\documentclass[10pt,a4paper]{amsart}
\usepackage{amsmath}
\usepackage[utf8]{inputenc}
\usepackage{amssymb,amsmath}
\usepackage{amsthm}
\usepackage{epsfig}
\usepackage{graphicx}
\usepackage{bm}
\usepackage{color}
\usepackage{bm}
\usepackage{amsrefs}
\usepackage{pdfpages}
\usepackage{lineno}
\usepackage{setspace} 
\usepackage{mathrsfs}
\input xypic
\usepackage{xcolor}
\usepackage{soul}
\usepackage{todonotes}
\setuptodonotes{inline}
\usepackage{mathtools}
\input xypic
\usepackage{vmargin}
\usepackage{ulem}
\setpapersize{A4}
\setmargins{2.5cm}       
{1.5cm}                        
{16.5cm}                      
{23.42cm}                    
{10pt}                           
{1cm}                           
{0pt}                             
{2cm}                           

\newtheorem{theorem}{Theorem}[section]
\newtheorem{definition}[theorem]{Definition}
\newtheorem{lemma}[theorem]{Lemma}
\newtheorem{corollary}[theorem]{Corollary}
\newtheorem{proposition}[theorem]{Proposition}
\newtheorem*{stheorem}{Theorem}
\theoremstyle{definition}
\newtheorem{remark}{Remark}[section]

\title{Central limit theorems for Soft random simplicial complexes }
\author{Juli\'an David Candela Coca}
\date{}
\thanks{This work was supported by the the CONAHCYT postgraduate studies scholarship number 940355 and by the grant N62909-19-1-2134 from the US Office of Naval Research Global and the Southern Office of Aerospace Research and Development of the US Air Force of Scientific Research.}
\begin{document}
\vspace{0.5cm}
\large
\begin{abstract}
\large
A soft random graph $G(n,r,p)$ can be obtained from the random geometric graph $G(n,r)$ by keeping every edge in $G(n,r)$ with probability $p$ \cite{MR3476631}. The soft random simplicial complexes is a model for random simplicial complexes built over the soft random graph $G(n,r,p)$. Similarly to \cite{MR3509567}, this new model depends on  a probability vector $\rho$ which allows the simplicial complexes to present randomness in all dimensions. In this article, we use a normal approximation theorem to prove  central limit theorems for the number of $k$-faces and for the  Euler's characteristic for soft random simplicial complexes.
\end{abstract}
\maketitle
\section{Introduction}

In probability theory, the central limit theorem (abbreviated  CLT) establishes that, for independent, identically distributed random variables, not necessarily normal, the  distribution of the normalized sample mean tends towards the standard normal distribution. One consequence of this result is that it implies that probabilistic and statistical methods that work for normal distributions can be used to address problems involving the sample mean of other types of distributions. A number of  authors have studied central limit theorems for different random combinatorial and topological models \cite{MR1986198},\cite{MR3079211},\cite{MR4252193}. One of them, the random geometric graph $G(n,r)$ \cite{MR1986198}, is a random graph over a set of $n$ uniformly distributed points on the set $\mathbb{R}^d$, in which each edge is included if the distance between its vertices is less than $r$. In \cite{MR1986198}, Penrose extensively studied the properties of this graph and  proved  central limit theorems for the number of components and for the number of subgraphs of $G(n,r)$ isomorphic to a given graph $\Gamma$.

Other models, such as the random Vietoris-Rips complex $R(n,r)$ and the random \v{C}ech complex $C(n,r)$, can be viewed as  higher dimensional generalizations of the random geometric graph $G(n,r)$. In \cite{MR3079211}, Kahle and Meckes studied the  limiting distributions of the Betti numbers of these two random geometric simplicial complexes by  proving central limit theorems for their Betti numbers.  In \cite{MR3509567}, Costa and Farber described a model for random simplicial complexes, built on a set on $n$ points, that depend on a probability vector  $\rho=(p_0,...,p_r)$ which allows the model to have randomness in every possible dimension. Costa and Farber, in their series of papers \cite{MR3509567},\cite{MR3661651},\cite{MR3604492}, studied important topological properties such as  connectivity, the fundamental group, and the behavior of the Betti numbers for this model.
More recently, in \cite{MR4252193}, the authors extended this research by proving limit theorems for topological invariants in  a dynamic version of Costa and Farber's model.

The soft random simplicial complexes $R(n,r,\rho)$ and $C(n,r,\rho)$, which we study in this paper, can be thought of as a combination of Costa and Farber's model and the random geometric complexes $R(n,r)$ and $C(n,r)$, or, equivalently, they can be seen as simplicial complexes built over the soft random graph $G(n,r,p)$ \cite{MR3476631} or  by taking $q=1-p_1$ and $p=0$,  over the $(p,q)$-perturbed random geometric graph $G_n^{q,p}$ \cite{MR4042102},\cite{MR4605137}. In this model, each simplicial complex is obtained from the original one by  keeping each $k$-face with probability $p_k$.  In \cite{yo}, the present author studied the expected topological properties  for this model in different regimens. We extend this work by proving  central limit theorems for  the number of $k$-faces and the Euler characteristic of the soft random simplicial complexes $R(n,r,\rho)$ and $C(n,r,\rho)$. 

The paper is organized as follows. Section\ref{section2} introduce notation, define the soft random simplicial complexes $R(n,r,\rho)$ and $C(n,r,\rho)$, and present the two  main results. In Section \ref{section3}, we address preliminary results as the asymptotic behavior of the variance and other important results that will be needed for calculations. The last two sections will be used to prove the main theorems:  The central limit theorem for the number of $k$-faces $f_k$ is proved in Section~\ref{section4}, while the Central limit theorem  for the Euler characteristic is proved in section \ref{section5}.

The results contained in this article form a part of my Ph.D. research at the Centro de Investigaci\'on en Matem\'aticas, A.C. under the supervision of Antonio Rieser.

\section{Notation and  Main Results}\label{section2}
The random geometric graph  $G(n,r)$ was  extensively studied by Penrose in his book Random geometric graphs \cite{MR1986198}. 

Given $n$ independent, identically distributed random variables  $X_1,\dots,X_n$ in $\mathbb{R}^d$ having a common probability density function $f:\mathbb{R}^d\rightarrow [0,\infty)$, the point process $\mathcal{X}_n$ is the union $\mathcal{X}_n=\bigcup\limits_{i=1}^n\{X_i\}$. Through this document, we assume that $f$  satisfies
$$\int_{\mathbb{R}^d}f(x)dx=1$$
and that for every $A\subset \mathbb{R}^d$ we have 
$$P(X_1\in A)=\int_{A}f(x)dx.$$
The random geometric graph $G(n,r)$, also denoted by $G(\mathcal{X}_n, r)$ is the the random  graph with vertex set $\mathcal{X}_n=\{X_1,\dots,X_n\}$ and an edge between $X_i$ and $X_j$  when $d(X_i, X_j)\leq r$, for $r=r(n)$ a function tending to zero. This graph was thoroughly studied by Penrose in \cite{MR1986198}.
On the other hand, we state the following coupling presented in
\cite{MR1986198}.  Given $\lambda>0$, let $N_\lambda$ be a Poisson
random variable independent of $\mathcal{X}_n=\{X_1,\dots,X_n\}$ and let
$$\mathcal{P}_\lambda= \{X_1,\dots,X_{N_\lambda}\}.$$
The following result from \cite{MR1986198} guarantees
that $\mathcal{P}_\lambda$ is in fact a Poisson process.
\begin{proposition}[Propositions 1.5 in \cite{MR1986198}]
$\mathcal{P}_\lambda$ is a Poisson process on $\mathbb{R}^d$ with intensity $\lambda f$.
\end{proposition}

This Poisson process is introduced to follow a common argument used by Penrose in \cite{MR1986198} in which he prove a result for the Poisson version of the random geometric graph $G(\mathcal{P}_n, r)$, that is the geometric graph with vertex set $\mathcal{P}_\lambda$ with $\lambda=n$  , and then using a de-Poissonization argument to conclude the same result for the random geometric graph $G(n,r)$.\\
$\mathcal{P}_\lambda$ is a Poisson process coupled with $\mathcal{X}_n$ related by the following result. 

\begin{theorem}[Theorem 1.6 in \cite{MR1986198}]\label{Penrose-5}
Let $\lambda>0$. Suppose $j\in \mathbb{M}$ and suppose $h(Y,X)$ is a bounded measurable function defined on all pairs of the form $(Y,X)$ with $X$ a finite subset of $\mathbb{R}^d$ and $Y$ a subset of $X$ satisfying $h(Y,X)=0$ except when $Y$ has $j$ elements. Then
$$E(\sum_{Y\subset \mathcal{P}_\lambda}h(Y,X))=\frac{\lambda^j}{j!}E(h(\mathcal{X'}_j, \mathcal{X'}_j\cup \mathcal{P}_\lambda))$$
where $\mathcal{X'}_j$ is an independent copy of $\mathcal{X}_j$, independent of $\mathcal{P}_\lambda$.
\end{theorem}

\begin{theorem}[Theorem 1.7 in \cite{MR1986198}]\label{Penrose-p}
Let $\lambda>0$. Suppose that $j_1,j_2\in \mathbb{N}$ and, for $i=1,2$, that $h_i(Y)$ is a bounded measurable function defined in all finite subset $Y \subset\mathbb{R}^d$ and satisfying $h_i(Y)=0$ except when $Y$ has $j_i$ elements. Then
$$E\left(\sum_{Y_1\subset \mathcal{P}_\lambda}\sum_{Y_2\subset \mathcal{P}_\lambda}h_1(Y_1)h_2(Y_2)1_{\{Y_\cap Y_2=\emptyset\}}\right)=E\left( \sum_{Y_1\subset \mathcal{P}_\lambda}h_1(Y_1)\right)E\left(\sum_{Y_2\subset \mathcal{P}_\lambda}h_2(Y_2)\right)$$
\end{theorem}

We now define some notation from \cite{MR1986198} that
we will use along this document. Given a finite point set
$Y\subset \mathbb{R}^d$, the first element of $Y$ according to the
lexicographic order will be called the left most point of $Y$ and will
be denoted by $LMP(Y)$. For an open set $A\subset \mathbb{R}^d$  with Lebesgue measure of its boundary $Leb(\partial A)=0$  and a connected graph $\Gamma$ on $k$ vertices,  let $G_{n,A}(\Gamma)$ and  $G'_{n,A}(\Gamma)$ be the random variables  that count the number of induced subgraphs of $G(n,r)$ isomorphic to $\Gamma$ and  the number of induced subgraphs of $G(\mathcal{P}_n,r)$ isomorphic to $\Gamma$, respectively for which the left most point of the vertex set lies in $A$. Similarly, for $\Gamma$ as above,  let $J_{n,A}(\Gamma)$ and $J'_{n,A}(\Gamma)$ be the random variables  that count the number of components of $G(n,r)$ isomorphic to $\Gamma$ and  the number components of $G(\mathcal{P}_n,r)$ isomorphic to $\Gamma$, respectively, for which the left most point of the vertex set lies in $A$.

On the other hand, given a
  connected graph $\Gamma$  of order $k$ and
sets $Y$ and $A$ as above, define the 
characteristic functions $h_\Gamma$ as:
\begin{equation}
h_{\Gamma}(Y)=1_{\{G(Y,1)\cong \Gamma\}},~~~h_{\Gamma,n,A}(Y)=1_{\{G(Y,r)\cong \Gamma\}\cap \{LMP(Y)\in A\}}
\end{equation}
and $\mu_{\Gamma,A}$ as:
\begin{equation}\label{mu}
\mu_{\Gamma, A}=\frac{1}{k!}\int_A f(x)^k~dx \int_{(\mathbb{R}^d)^{k-1}}h_{\Gamma}(\{0,x_1,\dots,x_{k-1}\})~dx_1dx_2\dots dx_{k-1}
\end{equation}
notice that $\mu_{\Gamma, A}$ is finite since it is bounded by a constant times $||f||Leb(B(0,k))$. In case $A=\mathbb{R}^d$, the  subscript $A$ in all  the notation stated above will be omitted. Also, when $\Gamma=K_m$ is the complete graph on $m$ vertices, we write $\mu_{m, A}$ for $\mu_{\Gamma, A}=\mu_{K_m, A}$.

We will now state some important results from \cite{MR1986198}:    

\begin{definition}
Suppose that $\Gamma$ is a graph of order $k\geq 2$. We say that $\Gamma$ is feasible iff $P(G(\mathcal{X}_k,r)\cong \Gamma)>0$ for some $r>0$.
\end{definition}

\begin{proposition}[Propositions 3.1 in \cite{MR1986198}]\label{Penrose-1}
Suppose that $\Gamma$ is a feasible connected graph of order $k\geq 2$ and $A\subset \mathbb{R}^d$ is open with $Leb(\partial A)=0$. Also, suppose that $r=r(n)\rightarrow 0$. Then
$$\lim_{n\rightarrow \infty}\frac{G_{n,A}(\Gamma)}{n^k r^{d(k-1)}}=\lim_{n\rightarrow \infty}\frac{G_{n,A}'(\Gamma)}{n^k r^{d(k-1)}}=\mu_{\Gamma,A},$$
where $\mu_{\Gamma,A}$ is the constant defined in (\ref{mu}),  which depends on $\Gamma$ and $k$.
\end{proposition}

\begin{proposition}[Propositions 3.2 in \cite{MR1986198}]\label{Penrose-2}
Let $\Gamma$ be a feasible connected graph of order $k\geq 2$ and let $J_{n,A}(\Gamma)$ be the number of components of $G(n,r)$ isomorphic to $\Gamma$ with LMP of the set of vertices in $A$. Also, suppose that $nr^d\rightarrow 0$. Then
$$\lim_{n\rightarrow \infty}\frac{J_{n,A}(\Gamma)}{n^k r^{d(k-1)}}=\lim_{n\rightarrow \infty}\frac{J'_{n,A}(\Gamma)}{n^k r^{d(k-1)}}=\mu_{\Gamma,A}$$
where $\mu_\Gamma$ a constant that depends on $\Gamma$ and $k$.
\end{proposition}
\begin{proposition}[Propositions 3.3 in \cite{MR1986198}]\label{Penrose-3}
Suppose that $\Gamma$ is a feasible connected graph of order $k\geq 2$. Also, suppose that $nr^d\rightarrow \lambda\in(0,\infty)$. Then
$$\lim_{n\rightarrow \infty}\frac{J_{n,A}(\Gamma)}{n}=\lim_{n\rightarrow \infty}\frac{J'_{n,A}(\Gamma)}{n}=c_{\Gamma, A}$$
where $c_{\Gamma, A}$ a constant that depends on $\Gamma$ and $A$.
\end{proposition}
Furthermore, it is noticed by Kahle in \cite{MR2770552}, that Proposition~\ref{Penrose-3} can be extended to work for $G_{n,A}(\Gamma )$ as follows:
\begin{proposition}\label{Kahle-1}
Suppose that $\Gamma$ is a feasible connected graph of order $k\geq 2$. Also, suppose that $nr^d\rightarrow \lambda\in(0,\infty)$. Then
$$\lim_{n\rightarrow \infty}\frac{G_{n,A}(\Gamma)}{n}=\lim_{n\rightarrow \infty}\frac{G'_{n,A}(\Gamma)}{n}=d_{\Gamma, A}$$
for $c_{\Gamma, A}$ a constant that depends on $\Gamma$ and $A$.
\end{proposition}

\begin{definition}
The random Vietoris-Rips complex $R(n,r)$, also written $R(\mathcal{X}_n,r)$, is the simplicial complex with vertex set $\mathcal{X}_n$, and a $k$-simplex $\sigma=(X_0,\dots,X_{k})\in C(n,r)$ iff
$$B(X_i,\frac{r}{2})\cap B(X_j,\frac{r}{2})\neq \emptyset$$
for every pair $X_i,X_j\in \sigma.$ Also, when  vertex set is $\mathcal{P}_n$, the Poisson version of the random Vietoris-Rips complex will be represented by $R(\mathcal{P}_n,r)$.
\end{definition}

\begin{definition}
The random \v{C}ech complex $C(n,r)$, also written
  $C(\mathcal{X}_n,r)$, is the simplicial complex with vertex set
  $\mathcal{X}_n$, and a  $k$-simplex $\sigma=(X_0,\dots,X_{k})\in C(n,r)$ iff
$$\bigcap_{X_i\in \sigma} B(X_i, \frac{r}{2})\neq \emptyset.$$
Also, when  vertex set is $\mathcal{P}_n$, the the Poisson version of the random \v{C}ech complex will be represented by $C(\mathcal{P}_n,r)$.
\end{definition}

Along this article, we will introduce random variables associated to different versions of random simplicial complexes.  Thus, we need some notation that allow us to distinguish whether we are talking about a random variable related to the Poisson version of a random simplicial complex or not; or if we are talking about a random variable related to a random simplicial complex that depends on a probability vector $\rho$.
Therefore, we make the following conventions: All random variables related to a  random simplicial complex that depends on a probability vector $\rho(n)=(p_1(n),p_2(n),...)$ will have the explicit dependency of this parameter. Also, all random variables related to a random graph or random simplicial complexes built over the  Poisson process $\mathcal{P}_n$ will have added the super index $\mathcal{P}$. Thus:

\begin{enumerate}
\item Random variables without an explicit dependency of the parameter $\rho$ and without the super index $\mathcal{P}$ are associated to the graph $G(n,r)$ or the simplicial complexes $R(n,r)$ and $C(n,r)$. In any case, to random graphs or random complexes that do not depend on a probability vector $\rho$ and whose vertex set are $n$ independent distributed points  $X_1,...,X_n$.
\item Random variables without an explicit dependency of the parameter $\rho$ and with the super index $\mathcal{P}$ are associated to the graph $G(\mathcal{P}_n,r)$ or the simplicial complexes $R(\mathcal{P}_n,r)$ and $C(\mathcal{P}_n,r)$. In any case, to random graphs or random complexes that do not depend on a probability vector $\rho$ and whose vertex set depends on a Poisson process  $\mathcal{P}_n$.
\item Random variables with an explicit dependency of the parameter $\rho$ and without the super index $\mathcal{P}$ are associated to $R(n,r,\rho)$ or $C(n,r,\rho)$. In any case, to  random complexes that do depend on a probability vector $\rho$ and whose vertex set are $n$ independent distributed points  $X_1,...,X_n$.
\item Random variables with an explicit dependency of the parameter $\rho$ and with the super index $\mathcal{P}$ are associated to $R(\mathcal{P}_n,r,\rho)$ or $C(\mathcal{P}_n,r,\rho)$. In any case, to  random complexes that do depend on a probability vector $\rho$ and whose vertex set depend on a  Poisson process $\mathcal{P}_n$.
\item In order to not overload the writing, we often omit the explicit dependence of the multiparameter $\rho$ on $n$.
\end{enumerate}

\begin{definition}\label{h's}
Consider $\mathcal{X}_n$ and $\mathcal{P}_n$ as before. Let $A\subset \mathbb{R}^d$ open, $\sigma$ a simplicial complex and  $\Gamma$ a graph. Let $\mathcal{H}\in\{\mathcal{X}_n,\mathcal{P}_n\}$ a set of vertices and $X\subset \mathcal{H}$ finite, define the following   characteristic functions:
\begin{enumerate}
\item $h_{\mathcal{H},r,\Gamma, A}(X)$  as the characteristic function that is 1 when $G(X, r)$  is a subgraph of $G(\mathcal{H}, r)$ isomorphic to $\Gamma$ with $LMP(X) \in A$.
\item $h_{\mathcal{H},r,\Gamma, A,p}(X)$ as the characteristic function that is 1 when $G(X, r, p)$ is a subgraph of $G(\mathcal{H}, r, p)$  isomorphic to $\Gamma$ with $LMP(X) \in A$. 
\item $h_{\mathcal{H},r,\sigma, A}(X)$ as the characteristic function that is 1 when $R(X, r)$  is a subcomplex of $R(\mathcal{H}, r)$  isomorphic to $\sigma$ with the $LMP(X) \in A$. 
\item $h_{\mathcal{H},r,\sigma,A,\rho}(X)$  as the characteristic function that is 1 when $R(X,r,\rho)$is a subcomplex of $R(\mathcal{X}_n,r, \rho)$  isomorphic to $\Gamma$ with the $LMP(X) \in A$. 

\item $\mathfrak{h}_{\mathcal{H},r,\sigma, A}(X)$ as the characteristic function that is 1 when $C(X, r)$ is a subcomplex of $C(\mathcal{X}_n, r)$  isomorphic to $\sigma$ with the $LMP(X) \in A$. 
\item $\mathfrak{h}_{\mathcal{H},r,\sigma,A,\rho}(X)$ as the characteristic function that is 1 when $C(X,r,\rho)$  is a subcomplex of $C(\mathcal{H},r, \rho)$  isomorphic to $\sigma$ with the the $LMP(X) \in A$.  
\end{enumerate}
\begin{remark}
 We also use the fallowing conventions in all notation above.
 \begin{itemize}
 \item If $A=\mathbb{R}^d$, we omit the subscript $A$.
\item If $\mathcal{H}=\mathcal{X}_n$, we replace the subscript $\mathcal{X}_n$ by $n$. 
\item  If the case that $\Gamma_{k+1}$ is the complete graph on $k+1$ vertices or $\sigma_{k+1}$ is a $k$-face, we replace the subscripts $\Gamma_{k+1}$ and $\sigma_{k+1}$ for $k+1$. 
 \end{itemize}   
\end{remark} 
\end{definition}

\begin{definition}
Let $f(n)$ and $g(n)$ real valued functions we say that they are asymptotically equal and write $f(g)\sim g(n)$ iff
$$\lim_{n\to \infty}\frac{f(n)}{g(n)}=1.$$
\end{definition}
\subsection{Main results}
We now introduce the random simplicial complexes and its principal results which we concern in these article.
\begin{definition}
Let $\rho(n)$ be the infinite multiparameter vector $\rho(n)=(p_1(n),p_2(n),\dots)$ with $0\leq p_i(n)\leq 1$ for each $i$. let $r=r(n)$ be a function tending to zero as $n$ goes to $\infty$. 
\begin{enumerate}
\item The soft random Vietoris-Rips complex $R(n,r,\rho)$ and its Poisson version $R(\mathcal{P}_n,r,\rho)$   are the simplicial complexes with vertex set $\mathcal{X}_n$ and $\mathcal{P}_n$, respectively,  in which a $k$-simplex $\sigma=(X_0,\dots,X_{k})$ in included in it with probability $p_k$ iff:
$$B(X_i,\frac{r}{2})\cap B(X_j,\frac{r}{2})\neq \emptyset$$
for every pair $X_i,X_j\in \sigma$. 
\item The soft random \v{C}ech complexes $C(n,r,\rho)$ and its Poisson version $C(\mathcal{P}_n,r,\rho)$   are the simplicial complexes with vertex set $\mathcal{X}_n$ and $\mathcal{P}_n$, respectively, in which a $k$-simplex $\sigma=(X_0,\dots,X_{k})$ in included on it with probability $p_k$ iff:
$$\bigcap_{X_i\in \sigma} B(X_i, \frac{r}{2})\neq \emptyset.$$
\end{enumerate}
\end{definition}

We present now the main results we want to address by this article, the central limit theorems for the number of k-faces $f_k$ and the Euler characteristic $\chi$ of the soft random simplicial complexes. 

\begin{stheorem}[\ref{Candela15}]
Let $R(n,r,\rho)$ be the soft random Vietoris-Rips complex and $C(n,r,\rho)$ the random soft \v{C}ech complex with multiparameter $\rho=(p_1,p_2,...,p_{n-1},...)$. Let $f_k(\rho)$ the number of $k$-faces of any of these two models. Then if $nr^d\to 0$  and $\prod\limits_{i=1}^k p_i^{\binom{k+1}{i+1}}n^{k+1}r^{dk}\to \infty$:
$$\frac{ f_k(\rho, R(n,r,\rho))-E(f_k(\rho, R(n,r,\rho)))}{\sqrt{\mu_{k+1}\prod\limits_{i=1}^kp_i^{\binom{k+1}{i+1}}n^{k+1}r^{dk}}}\xRightarrow{D} N(0, 1)$$
and
$$\frac{ f_k(\rho, C(n,r,\rho))-E(f_k(\rho, C(n,r,\rho)))}{\sqrt{\nu_{k+1}\prod\limits_{i=1}^kp_i^{\binom{k+1}{i+1}}n^{k+1}r^{dk}}}\xRightarrow{D} N(0, 1)$$
\end{stheorem}

\begin{stheorem}[\ref{Candela16}]
Let $\Gamma$ one the soft random simplicial complexes $\Gamma\in\{R(n,r,\rho), C(n,r,\rho)\}$ with multiparameter $\rho=(p_1,p_2,...,p_{n-1},...)$ and let $\chi(\rho,\Gamma)$ be its Euler characteristic. Assume that $nr^d\to 0$ and that there exists a non-negative integer $l$ such that 
$$\prod\limits_{i=1}^{l+1} p_i^{\binom{l+2}{i+1}}n^{l+2}r^{d(l+1)}\to 0~~  and~~\prod\limits_{i=1}^l p_i^{\binom{l+1}{i+1}}n^{l+1}r^{dl}\to \infty$$
then 
$$ \frac{\chi(\rho,\Gamma)-E(\chi(\rho,\Gamma))}{\sqrt{n}}\xRightarrow{D} N(0,1).$$
\end{stheorem}

\begin{remark}
\begin{enumerate}
\item[]
\item If, for each $i$, $p_i=1$ then the hypothesis of Theorem~\ref{Candela16} can be achieved by taking $r^d\leq \frac{1}{n^a}$ for $a\in (\frac{l+2}{l+1}, \frac{l+1}{l})$ 
\item Let $\epsilon>0$ given. If $\prod\limits_{i=1}^l p_i^{\binom{l+1}{i+1}}<\frac{1}{n^b}$, then the hypothesis of Theorem~\ref{Candela16} can be achieved by taking $r^d\leq \frac{1}{n^a}$ for $a=(1+\epsilon)$ and $b\in(1-\epsilon(l+1), 1-\epsilon l)$. 
\end{enumerate}
\end{remark}

\section{Preliminary Results}\label{section3}
In this section, we present a number of preliminary results which we will use  to prove the two main theorems of this article. The reader may skip this section on a first reading and refer back to it as needed. 
\begin{theorem}\label{equalh} Let $A\subset \mathbb{R}^d$ be open, $k\geq 1$,  $\sigma$ a $k$-simplex and denote by $\Gamma_{k+1}$ the complete graph of $k+1$ vertices. For $h_{n,r,\Gamma_{k+1},A,}$ and $h_{n,r,\sigma,A}$ as in Definition~\ref{h's}, we have, for any vertex set
$\mathcal{X}_{k+1}=\{X_1,\dots ,X_{k+1}\}$
$$h_{n,r,\sigma,A}(\mathcal{X}_{k+1}) = h_{n,r,\Gamma_{k+1},A}(\mathcal{X}_{k+1}),$$
and therefore
$$E(h_{n,r,\sigma,A}(\mathcal{X}_{k+1}))=E(h_{n,r,\Gamma_{k+1},A}(\mathcal{X}_{k+1}))$$
\end{theorem}
\begin{proof}
 $\mathcal{X}_{k+1}$ spans a $k$-face $\sigma\in R(n,r)$ iff $\mathcal{X}_{k+1}$ spans the complete graph  $\Gamma_{k+1}$ in $G(n,r)$. Then, $h_{n,r,\sigma,A}(\mathcal{X}_{k+1}) = h_{n,r,\Gamma_{k+1},A}(\mathcal{X}_{k+1})$. 
\end{proof}
The following constants will appear very often during the proofs. In particular, they are  constant associated to the expected value of the number of $k$-faces in the simplicial complexes $R(n,r)$ and $C(n,r).$

\begin{definition}\label{numudef}
Let $A\subset \mathbb{R}^d$ be open and $k\geq1$. Consider the notation of Definition~\ref{h's}.
\begin{align}
\mu_{0, A}=&\nu_{0,A}=\int_{A}f(x)dx.\\
\mu_{k+1, A}=&\frac{1}{(k+1)!}\int_A f(x)^kdx \int_{(\mathbb{R}^d)^{k}}h_{n,1,k+1}(\{0,x_1,...,x_{k}\})~dx_1dx_2...dx_{k}.\label{newmu}\\
\nu_{k+1, A}=&\frac{1}{(k+1)!}\int_A f(x)^kdx \int_{(\mathbb{R}^d)^{k}}\mathfrak{h}_{n,1,k+1}(\{0,x_1,...,x_{k}\})~dx_1dx_2...dx_{k}.\label{newnu}
\end{align}
When $A=\mathbb{R}^d$, we omit the subscript $A$. 
\end{definition}

\begin{theorem}\label{newcosito}
For $k\geq 1$, let $X$ a  $k+1$ simplex.  If $G_{n,A}(X)$ and $J_{n,A}(X)$ are, respectively, the number of subcomplex and components of $C(n,r)$ isomorphic to $X$, with left most point of the $k$-simplex in $A$, then  
$$\lim_{n\to \infty}\frac{E(G_{n,A}(X))}{n^{k+1}r^{d(k)}}=\lim_{n\to \infty}\frac{E(J_{n,A}(X))}{n^{k+1}r^{d(k1)}}=\nu_{k+1, A}.$$
\end{theorem}
\begin{proof}
\begin{align*}
E(G_{n,A}(X))&=\binom{n}{k+1}\int_{\mathbb{R^d}}...\int_{\mathbb{R^d}}\mathfrak{h}_{n,r,A}(x_1,...,x_{k+1})\prod_{i=1}^{k+1}f(x_i)dx_1...dx_{k+1}\\
&=\int_{\mathbb{R^d}}...\int_{\mathbb{R^d}}\mathfrak{h}_{n,r,A}(x_1,...,x_{k+1})f(x_1)^{k+1} dx_1...dx_{k+1}\\
 &+\int_{\mathbb{R^d}}...\int_{\mathbb{R^d}}\mathfrak{h}_{n,r,A}(x_1,...,x_{k+1})\left(\prod_{i=1}^{k+1}f(x_i)-f(x_1)^{k}\right)dx_1...dx_{k+1}.   
\end{align*}
Now, as in proof of Proposition 1.3 from \cite{MR1986198},
$$ \int_{\mathbb{R^d}}...\int_{\mathbb{R^d}}\mathfrak{h}_{n,r,A}(x_1,...,x_{k+1})\left(\prod_{i=1}^{k+1}f(x_i)-f(x_1)^{k+1}\right)dx_1...dx_{k+1}\to 0.$$
For the first part, making the change of variables $x_i=x_1+ry_i$ for $2\leq i\leq k+1$:
\begin{align*}
\int_{\mathbb{R^d}}...\int_{\mathbb{R^d}}&\mathfrak{h}_{n,r,A}(x_1,...,x_{k+1})f(x_1)^{k+1} dx_1...dx_{k+1}\\
&=r^{d(k)}\int_{\mathbb{R^d}}...\int_{\mathbb{R^d}}\mathfrak{h}_{n,r,A}(x_1,x_1+ry_2,...,x_1+ry_{k+1})f(x_1)^{k+1} dy_2...dy_{k+1}dx_1
\end{align*}
again, as in the proof of Proposition 1.3 from \cite{MR1986198}, the above is asymptotic to $(k+1)!\nu_{k+1, A}r^{dk}$. Thus 
$$E(G_{n,r}(X))\sim \binom{n}{k+1}\nu_{X}r^{dk} \sim \frac{n^k}{k!}\nu_{k+1, A}r^{dk}.$$
On the other hand, let $\hat{\mathfrak{h}}_{n,r,A}(x_1,\dots,x_{k+1})$ the characteristic function that is $1$ if the set $\{x_1,\dots,x_{k+1}\}$ spans a component of $C(n,r)$ isomorphic to $X$. Then
$$E(J_{n,r}(X))=\binom{n}{k+1}P(\hat{\mathfrak{h}}_{n,r,A}(x_1,...,x_{k+1})=1).$$
Now, the conditional probability $Q$ of the set $\{x_1,...,x_{k+1}\}$ spanning a component isomorphic to $X$ in $C(n,r)$ given that spans a subcomplex isomorphic to $X$ satisfies that
$$Q=P(\hat{\mathfrak{h}}_{n,r,A}(x_1,...,x_{k+1})=1\mid \mathfrak{h}_{n,r,A}(x_1,...,x_{k+1})=1 )\to 1$$
since 
$$Q\geq P(\text{no vertex $x_j \notin \{x_1,...,x_{k+1}\}$ is connected with $\{x_1,...,x_{k+1}\}$} ).$$
On the other hand, the probability $Q_j$ that a vertex $x_j \notin \{x_1,...,x_{k+1}\}$ is connected to a vertex in $\{x_1,...,x_{k+1}\}$ satisfies
$$Q_j=\int_{\bigcup\limits_{i=1}^{k+1}B(x_i,r)}f(x_j)dx_j\leq ||f||_{\infty}\int_{B(x_1,(k+1)r)}dx_j \leq ||f||_{\infty}((k+1)r)^d\int_{B(0,1)}dx_j $$
since $X$ is connected. Therefore, if $\theta_d=Leb(B(0,1))$, then by independence 
$$Q\geq (1-Q_j)^{n-k-1}\geq (1-||f||_{\infty}((k+1)r)^d \theta_d)^{n-k-1}\to 1.$$
Thus
\begin{align*}
E(J_{n,A}(X))&=\binom{n}{k+1}P(\hat{\mathfrak{h}}_{n,r,A}(x_1,...,x_{k+1})=1)\\
&=\binom{n}{k+1}P(\mathfrak{h}_{n,r,A}(x_1,...,x_{k+1})=1)Q\\
&=\binom{n}{k+1}P(\mathfrak{h}_{n,r,A}(x_1,...,x_{k+1})=1)\\
&= E(G_{n,A}(X))
\end{align*}

\end{proof}

\begin{theorem} Let $A\subset \mathbb{R}^d$ open and $k\geq 1$. For $h_{n,r,\sigma,A,\rho}$ and $h_{n,r,\sigma,A,\rho}$ as in Definition~\ref{h's}, there for a constant $c=\frac{\nu_{k+1, A}}{\mu_{k+1, A}}$:
$$E(\mathfrak{h}_{n,r,k+1,A}(\mathcal{X}_{k+1}))\sim cE(h_{n,r,k+1,A}(\mathcal{X}_{k+1}))).$$
\end{theorem}
\begin{proof} 
Since the vertices $\mathcal{X}_{k+1}=\{X_1,\dots ,X_{k+1}\}$  have a common probability density function $f:\mathbb{R}^d\rightarrow [0,\infty)$, then
\begin{align*}
  E(h_{n,r,\Gamma,A}(\mathcal{X}_{k+1}))&=\int_{\mathbb{R^d}}...\int_{\mathbb{R^d}}h_{n,r,k+1,A}(x_1,...,x_{k+1})\prod_{i=1}^{k+1}f(x_i)dx_1...dx_{k+1}\\
&=\int_{\mathbb{R^d}}...\int_{\mathbb{R^d}}h_{n,r,k+1,A}(x_1,...,x_{k+1})f(x_1)^{k+1} dx_1...dx_{k+1}\\
 &\quad+\int_{\mathbb{R^d}}...\int_{\mathbb{R^d}}h_{n,r,k+1,A}(x_1,...,x_{k+1})\left(\prod_{i=1}^{k+1}f(x_i)-f(x_1)^{k+1}\right)dx_1...dx_{k+1}.  
\end{align*}
We now observe that, as in the proof of Proposition 1.3 from \cite{MR1986198},
as $n\to \infty$
$$ \int_{\mathbb{R^d}}...\int_{\mathbb{R^d}}h_{n,r,\Gamma,A}(x_1,...,x_{k+1})\left(\prod_{i=1}^{k+1}f(x_i)-f(x_1)^{k+1}\right)dx_1...dx_{k+1}\to 0$$
since as $n\to \infty$, $r\to 0$ and therefore  $|\prod\limits_{i=1}^{k+1}f(x_i)-f(x_1)^{k+1}|\to 0$.
On the other hand, making the change of variables $x_i=x_1+ry_i$ for $2\leq i\leq k+1$:
\begin{align*}
&\int_{\mathbb{R^d}}...\int_{\mathbb{R^d}}h_{n,r,\Gamma,A}(x_1,...,x_{k+1})f(x_1)^{k+1} dx_1...dx_{k+1}&\\
&\qquad \qquad=r^{dk}\int_{\mathbb{R^d}}...\int_{\mathbb{R^d}}h_{n,r,\Gamma,A}(x_1,x_1+ry_2,...,x_1+ry_{k+1})f(x_1)^k dy_2...dy_{k+1}dx_1&    
\end{align*}
again, following the same steps as in the proof of Proposition 1.3 from \cite{MR1986198}, the above is asymptotic to $\mu_{k+1,A}r^{dk}$. Thus 
$$E(h_{n,r,k+1,A}(\mathcal{X}_{k+1}))\sim (k+1)!\mu_{k+1,A}r^{dk}.$$
Along  the same lines, we can prove that 
$$E(\mathfrak{h}_{n,r,k+1,A}(\mathcal{X}_{k+1}))\sim(k+1)!\nu_{k+1, A} r^{dk}.$$
Therefore, taking $c=\frac{\nu_{k+1, A}}{\mu_{k+1, A}}$,
$$\lim_{n\to \infty}\frac{E(h_{n,r,k+1,A}(\mathcal{X}_{k+1}))}{E(\mathfrak{h}_{n,r,k+1,A}(\mathcal{X}_{k+1}))}=1.$$
\end{proof}

\begin{theorem}\label{hvsh}
Let $A\subset \mathbb{R}^d$ be open and $k\geq1$.  Then
$$E(h_{n,r,k+1,A,\rho}(\mathcal{X}_{k+1}))=\prod_{i=1}^{l}p_i^{\binom{k+1}{i+1}} E(h_{n,r,k+1,A}(\mathcal{X}_{k+1}))$$ 
and
$$E(\mathfrak{h}_{n,r,k+1,A,\rho}(\mathcal{X}_{k+1}))=\prod_{i=1}^{l}p_i^{\binom{k+1}{i+1}} E(\mathfrak{h}_{n,r,k+1,A}(\mathcal{X}_{k+1})),$$
where $h_{n,r,k+1,A,\rho}(\mathcal{X}_{k+1})$ and $\mathfrak{h}_{n,r,k+1,A}(\mathcal{X}_{k+1})$  are defined in Definition~\ref{h's}.
\end{theorem}

\begin{proof}
Consider the conditional probability the set of points $\mathcal{X}_{k+1}$ spans a $k$-face in $R(n,r,\rho)$ given that it spans a $k$-face in $R(n,r)$. We observe that
$$P(h_{n,r,k+1,A,\rho}(\mathcal{X}_{k+1})=1\mid h_{n,r,k+1,A}(\mathcal{X}_{k+1})=1)= \prod_{i=1}^{k}p_i^{\binom{k+1}{i+1}}.$$
On the other hand, the event that $\mathcal{X}_{k+1}$ spans a $k$-face  in $R(n,r,\rho)$ is contained in the event that $\mathcal{X}_{k+1}$ spans a $k$-face in $R(n,r)$, thus:
$$P\left(h_{n,r,k+1,A,\rho}(\mathcal{X}_{k+1})=1\right)=P\left(h_{n,r,k+1,A}(\mathcal{X}_{k+1})=1\right) \prod_{i=1}^{k}p_i^{\binom{k+1}{i+1}}.$$
Therefore
$$E(h_{n,r,k+1,A,\rho}(\mathcal{X}_{k+1}))=E(h_{n,r,k+1,A}(\mathcal{X}_{k+1})\prod_{i=1}^{k}p_i^{\binom{k+1}{i+1}}.$$
The proof  for  $\mathfrak{h}_{n,r,k+1,A}$ is analogous.
\end{proof}

\begin{theorem}\label{prodh's}
Let $A\subset \mathbb{R}^d$ open and $k,l\geq 1$. Assume $l\geq k$ and that $\mathcal{X}_{k+1}$ and $\mathcal{X}_{l+1}$ have $j$ points in common. Then 
$$E(h_{n,r,k+1,A,\rho}(\mathcal{X}_{k+1})h_{n,r,l+1,A,\rho}(\mathcal{X}_{l+1}))=\prod_{i=1}^{l}p_i^{\binom{k+1}{i+1}+  \binom{l+1}{i+1}- \binom{j}{i+1}} E(h_{n,r,k+1,A}(\mathcal{X}_{k+1})h_{n,r,l+1,A}(\mathcal{X}_{l+1}))$$ 
and 
$$E(\mathfrak{h}_{n,r,k+1,A,\rho}(\mathcal{X}_{k+1})\mathfrak{h}_{n,r,l+1,A,\rho}(\mathcal{X}_{l+1}))
=\prod_{i=1}^{l}p_i^{\binom{k+1}{i+1}+  \binom{l+1}{i+1}- \binom{j}{i+1}} E(\mathfrak{h}_{n,r,k+1,A}(\mathcal{X}_{k+1})\mathfrak{h}_{n,r,l+1,A}(\mathcal{X}_{l+1})),$$
where $h_{n,r,k+1,A,\rho}(\mathcal{X}_{k+1})$ and $\mathfrak{h}_{n,r,k+1,A}(\mathcal{X}_{k+1})$  are defined in Definition~\ref{h's}.
\end{theorem}

\begin{proof}
The product $h_{n,r,k+1,A,\rho}(\mathcal{X}_{k+1})h_{n,r,l+1,A,\rho}(\mathcal{X}_{l+1})$ is 1 if and only if  $\mathcal{X}_{k+1}$ spans a $k$-face and $\mathcal{X}_{l+1}$ spans an  $l$-face in $R(n,r,\rho)$. On the other,
\begin{align*}
P(h_{n,r,k+1,A,\rho}(\mathcal{X}_{k+1})h_{n,r,l+1,A,\rho}(\mathcal{X}_{l+1})&=1\mid  h_{n,r,k+1,A}(\mathcal{X}_{k+1})h_{n,r,l+1,A}(\mathcal{X}_{l+1})=1)&\\
&= \prod_{i=1}^{l}p_i^{\binom{k+1}{i+1}+  \binom{l+1}{i+1}- \binom{j}{i+1}},
\end{align*}
and the event that $\mathcal{X}_{k+1}$ spans a $k$-face and $\mathcal{X}_{l+1}$ spans an  $l$-face in $R(n,r,\rho)$ is contained in the event that $\mathcal{X}_{k+1}$ spans a $k$-face and $\mathcal{X}_{l+1}$ spans an  $l$-face in $R(n,r)$. Therefore,
\begin{align*}
&P(h_{n,r,k+1,A,\rho}(\mathcal{X}_{k+1})h_{n,r,l+1,A,\rho}(\mathcal{X}_{l+1})=1)\\
&\qquad=P(h_{n,r,k+1,A}(\mathcal{X}_{k+1})h_{n,r,l+1,A}(\mathcal{X}_{l+1})=1) \prod_{i=1}^{l}p_i^{\binom{k+1}{i+1}+  \binom{l+1}{i+1}- \binom{j}{i+1}},
\end{align*}
and it follows that
\begin{align*}
E(h_{n,r,k+1,A,\rho}(\mathcal{X}_{k+1})h_{n,r,l+1,A,\rho}(\mathcal{X}_{l+1}))
=E(h_{n,r,k+1,A}(\mathcal{X}_{k+1})h_{n,r,l+1,A}(\mathcal{X}_{l+1})) \prod_{i=1}^{l}p_i^{\binom{k+1}{i+1}+  \binom{l+1}{i+1}- \binom{j}{i+1}}.
\end{align*}
The proof  for  $\mathfrak{h}$ is analogous.
\end{proof}

\begin{theorem}\label{prodhx1}
Let $A\subset \mathbb{R}^d$ be open. Assume $k\geq 1$. Then 
$$E(h_{n,r,k+1,A,\rho}(\mathcal{X}_{k}\cup\{x\})\cdot 1_{x \in A}(x))=\prod_{i=1}^{k}p_i^{\binom{k+1}{i+1}} E(h_{n,r,k+1,A}(\mathcal{X}_{k}\cup\{x\})\cdot 1_{x \in A})$$ 
and 
$$E(\mathfrak{h}_{n,r,k+1,A,\rho}(\mathcal{X}_{k}\cup\{x\})\cdot 1_{x \in A}(x))=\prod_{i=1}^{k}p_i^{\binom{k+1}{i+1}} E(\mathfrak{h}_{n,r,k+1,A}(\mathcal{X}_{l}\cup\{x\})\cdot 1_{x \in A})$$
where $1_{x \in A}(x)$ is the characteristic function of the set $A$ and  $h_{n,r,k+1,A,\rho}(\mathcal{X}_{k+1})$ and the functions $\mathfrak{h}_{n,r,k+1,A}(\mathcal{X}_{k+1})$  are defined in Definition~\ref{h's}.
\end{theorem}

\begin{proof}
The product $h_{n,r,k+1,A,\rho}(\mathcal{X}_{l}\cup\{x\})\cdot 1_{x \in A}$ is 1 if and only if $\mathcal{X}_{k}\cup\{x\}$ spans a $k$-face in $R(n,r,\rho)$ where the left most point of the set $\mathcal{ X}_{k}\cup\{x\}$  in $A$ and additionally $x\in A$. On the other hand, the conditional probability  that $\mathcal{X}_{k}\cup\{x\}$ spans a $k$-face in $R(n,r,\rho)$ for which its $LMP(\mathcal{ X}_{k}\cup\{x\})\in A$ and $x\in A$ given that $\mathcal{ X}_{k}\cup\{x\}$ spans a $k$-face in $R(n,r)$ for which its $LMP(\mathcal{ X}_{k}\cup\{x\})\in A$ and $x\in A$ is given by 

$$P(1_{\{x\in A\}}h_{n,r,k+1,A,\rho}(\mathcal{X}_{k}\cup\{x\})=1\mid 1_{\{x\in A\}}h_{n,r,k+1,A}(\mathcal{X}_{k}\cup\{x\})=1)= \prod_{i=1}^{l}p_i^{\binom{k+1}{i+1}}.$$
Furthermore, the event that $\mathcal{X}_{l}\cup\{x\}$ spans a $l$-face in $R(n,r,\rho)$ with $x\in A$ is contained in the event that $\mathcal{X}_{l}\cup\{x\}$ spans a $l$-face in $R(n,r)$ with $x\in A$, then:
$$P(1_{\{x\in A\}}h_{n,r,k+1,A,\rho}(\mathcal{X}_{l}\cup\{x\})=1)=P(1_{\{x\in A\}}h_{n,r,k+1,A}(\mathcal{X}_{l}\cup\{x\})=1) \prod_{i=1}^{l}p_i^{\binom{k+1}{i+1}},$$
and it follows that
$$E(1_{\{x\in A\}}h_{n,r,k+1,A,\rho}(\mathcal{X}_{l}\cup\{x\}))=\prod_{i=1}^{l}p_i^{\binom{l+1}{i+1}} E(1_{\{x\in A\}}h_{n,r,k+1,A}(\mathcal{X}_{l}\cup\{x\})).$$ 
The proof  for  $\mathfrak{h}$ is analogous.
\end{proof}
In the following result, we obtain an expression for the asymptotic behavior of the variance of the number of $k$-faces $f_{k,A}$, for each $k\geq 1$, of the simplicial complexes $R(n,r,\rho)$ and $R(\mathcal{P}_n,r,\rho)$ where left most point of the set of vertices (LMP) of the face is in $A$. For that we need the following notation:
\begin{definition} Let $A\subset \mathbb{R}^d$ be open and $k,l\geq 1$. Let $A\subset \mathbb{R}^d$ and  $(x_1,..., x_{k+l+2-j})\in \mathbb{R}^{d(k+l+2-j)}$. Using the conventions in Definition \ref{h's}, define the functions $h^j_{n,r,k+1,l+1,A,\rho}$ and $\mathfrak{h}^j_{n,r,k+1,l+1,A,\rho}$ as
$$h^j_{n,r,k+1,l+1,A,\rho}(x_1,...,x_{k+l+2-j})=h_{n,r,k+1,A,\rho}(x_1,...,x_{k+1}) h_{n,r,l+1,A,\rho}(x_1,...,x_j, x_{k+2},...,x_{k+l+2-j})$$
and 
$$\mathfrak{h}^j_{n,r,k+1,l+1,A,\rho}(x_1,...,x_{k+l+2-j})=\mathfrak{h}_{n,r,k+1,A,\rho}(x_1,...,x_{k+1}) \mathfrak{h}_{n,r,l+1,A,\rho}(x_1,...,x_j, x_{k+2},...,x_{k+l+2-j}).$$

\end{definition}
We now focus on proving several technical results that will allow us to determine the asymptotic behavior of the variance for the soft random simplicial complexes $R(n,r,\rho)$, $C(n,r,\rho)$, and their Poisson versions $R(\mathcal{P}_n,r,\rho)$ and $C(\mathcal{P}_n,r,\rho)$. In this first result, we consider the soft random Vietoris-Rips complex.
\begin{theorem}\label{Candela13}
Let $A\subset \mathbb{R}^d$ be an open set  and $k,l\geq 1$. Consider the soft random simplicial complexes  $R(n,r,\rho)$ and $C(n,r,\rho)$ with the multiparameter $\rho=(p_1,p_2,\dots)$ and the Poisson versions $R(\mathcal{P}_n,r,\rho)$ and $C(\mathcal{P}_n,r,\rho)$. Let $f_{k,A}(\rho)$ and $f_{k,A}^\mathcal{P}(\rho)$ be the number of $k$-faces of $R(n,r,\rho)$ and $R(\mathcal{P}_n,r,\rho)$, respectively, such that the leftmost point of the vertex set of the $k$-face is in $A$. Also, for convenience, we define
$$\Phi_{j,A}(f_k,f_l)=\frac{\int_{A}f(x)^{k+l+2-j}dx}{j!(k+1-j)!(l+1-j)!} \int_{\mathbb{R}^d}...\int_{\mathbb{R}^d}h^j_{n,1,k+1,l+1}(0,x_2,...,x_{k+l+2-j})dx_2\dots dx_{k+l+2-j}.$$
Then, if $nr^d\to 0$, we have
\begin{align*}
cov(f_{k,A}^\mathcal{P}(\rho),f_{l,A}^\mathcal{P}(\rho))=\sum_{j=1}^{\min(k+1,l+1)}n^{k+l+2-j}r^{d(k+l+1-j)}\prod_{i=1}^l p_i^{\binom{k+1}{i+1}+ \binom{l+1}{i+1}- \binom{j}{i+1}}\Phi_{j,A}(f_k,f_l)
\end{align*}
\end{theorem}

\begin{proof}
Without loss of generality, assume $k\leq l$, thus
\begin{align*}
cov(f_{k,A}^\mathcal{P}(\rho),f_{l,A}^\mathcal{P}(\rho))&=E(f_{k,A}^\mathcal{P}(\rho)f_{l,A}^\mathcal{P}(\rho))- E(f_{k,A}^\mathcal{P}(\rho))E(f_{l,A}^\mathcal{P}(\rho))\\
&=\sum_{j=0}^{k+1}E\left(\sum_{Y\subset \mathcal{P}_n}\sum_{Y'\subset \mathcal{P}_n} h_{\mathcal{P}_n,r,k+1,A,\rho}(X)h_{\mathcal{P}_n,r,l+1,A,\rho}(X')1_{\{|X\cap X'|=j\}}\right)\\
&\qquad- E(f_{k,A}^\mathcal{P}(\rho))E(f_{l,A}^\mathcal{P}(\rho))   
\end{align*}
with the functions $h_{\mathcal{P}_n,r,l+1,A,\rho}$ and $h_{\mathcal{P}_n,r,k+1,A,\rho}$ as in Definition~\ref{h's}. On the other hand,  by Theorem~\ref{Penrose-p}, the term $j=0$ on the left-hand side of the expression above corresponds to $E(f_{k,A}^\mathcal{P}(\rho))E(f_{l,A}^\mathcal{P}(\rho))$, which cancels the term on the right. We therefore have,
$$cov(f_{k,A}^\mathcal{P}(\rho),f_{l,A}^\mathcal{P}(\rho))=\sum_{j=1}^{k+1}E\left(\sum_{Y\subset \mathcal{P}_n}\sum_{Y'\subset \mathcal{P}_n} h_{\mathcal{P}_nn,r,k+1,A,\rho}(Y)h_{\mathcal{P}_n,r,l+1,A,\rho}(Y')1_{\{|Y\cap Y'|=j\}}\right)$$
thus
\begin{equation}\label{eq5.1}
cov(f_{k,A}^\mathcal{P}(\rho),f_{l,A}^\mathcal{P}(\rho))=\sum_{j=1}^{k+1}E\left(\sum_{X\subset \mathcal{P}_n} H(X)\right),
\end{equation}
where 
$$H_j(X)=1_{\{|X|=k+l+2-j\}}\sum_{Y\subset X}\sum_{Y'\subset X} h_{\mathcal{P}_n,r,k+1,A,\rho}(Y)h_{\mathcal{P}_n,r,l+1,A,\rho}(Y')1_{\{|Y\cap Y'|=j\}})$$
using Palm theory, Theorem~\ref{Penrose-5}, applied to Equation (\ref{eq5.1}), we have
$$E\left(\sum_{X\subset \mathcal{P}_n} H_j(X)\right)=\frac{n^{k+l+2-j}}{(k+l+2-j)!} E(H_j(\mathcal{X}_{k+l+2-j})).$$
Furthermore,
\begin{align*}
&E(H(\mathcal{X}_{k+l+2-j}))=E\left(\sum_{Y\subset \mathcal{X}_{k+l+2-j}}\sum_{Y'\subset \mathcal{X}_{k+l+2-j}} h_{n,r,k+1,A,\rho}(Y)h_{n,r,l+1,A,\rho}(Y')1_{\{|Y\cap Y'|=j\}})\right)\\
&=\binom{k+l+2-j}{k+1}\binom{k+1}{j}E(h_{n,r,k+1,A,\rho}(\mathcal{X}_{k+1})h_{n,r,l+1,A,\rho}(\mathcal{X}_{j}\cup (\mathcal{X}_{k+l+2-j}\backslash \mathcal{X}_{k+1})))\\
&=\frac{(k+l+2-j)!}{(k+1-j)!(l+1-j)!j!}E(h_{n,r,k+1,A,\rho}(\mathcal{X}_{k+1})h_{n,r,l+1,A,\rho}(\mathcal{X}_{j}\cup (\mathcal{X}_{k+l+2-j}\backslash \mathcal{X}_{k+1}))).
\end{align*}
Thus
\begin{align*}
cov(f_{k,A}^\mathcal{P}(\rho),f_{l,A}^\mathcal{P}(\rho))&=\sum_{j=1}^{k+1}E\left(\sum_{X\subset \mathcal{P}_n} H_j(X)\right)\\
&=\sum_{j=1}^{k+1}\frac{n^{k+l+2-j}}{(k+l+2-j)!}\frac{(k+l+2-j)!}{(k+1-j)!(l+1-j)!j!}\\
&\qquad\times E\left(h_{n,r,k+1,A,\rho}(\mathcal{X}_{k+1})h_{n,r,l+1,A,\rho}(\mathcal{X}_{j}\cup (\mathcal{X}_{k+l+2-j}\backslash \mathcal{X}_{k+1})))\right)\\
&=\sum_{j=1}^{k+1}\frac{n^{k+l+2-j}}{(k+1-j)!(l+1-j)!j!}\\
&\qquad\times E(h_{n,r,k+1,A,\rho}(\mathcal{X}_{k+1})h_{n,r,l+1,A,\rho}(\mathcal{X}_{j}\cup (\mathcal{X}_{k+l+2-j}\backslash \mathcal{X}_{k+1})))
\end{align*}
by Theorem~\ref{prodh's}:
\begin{align*}
&\qquad\qquad\qquad=\sum_{j=1}^{k+1}\frac{n^{k+l+2-j}}{(k+1-j)!(l+1-j)!j!}\prod_{i=1}^l p_i^{\binom{k+1}{i+1}+ \binom{l+1}{i+1}- \binom{j}{i+1}}\\
&\qquad\qquad\qquad\qquad\times E(h_{n,r,l+1,A}(\mathcal{X}_{k+1})h_{n,r,l+1,A}(\mathcal{X}_{j}\cup (\mathcal{X}_{k+l+2-j}\backslash \mathcal{X}_{k+1})))
\end{align*}
Therefore,
\begin{align}\label{sintag}
cov(f_{k,A}^\mathcal{P}(\rho),f_{l,A}^\mathcal{P}(\rho))&=\sum_{j=1}^{k+1}\frac{n^{k+l+2-j}}{(k+1-j)!(l+1-j)!j!}\prod_{i=1}^l p_i^{\binom{k+1}{i+1}+ \binom{l+1}{i+1}- \binom{j}{i+1}}\\
&\notag\quad\times E(h_{n,r,l+1,A}(\mathcal{X}_{k+1})h_{n,r,l+1,A}(\mathcal{X}_{j}\cup (\mathcal{X}_{k+l+2-j}\backslash \mathcal{X}_{k+1})))
\end{align}
On the other hand, It is proved in Proposition 3.7 in \cite{MR1986198} that for
\begin{align}\label{kik2}
c_r = r^{d(k+l+1-j)}\int_{A}f(x)^{k+l+2-j}dx \int_{\mathbb{R}^d}...\int_{\mathbb{R}^d}h^j_{n,1,k+1,l+1}(0,x_2,...,x_{k+l+2-j})dx_2\dots dx_{k+l+2-j}    
\end{align}
we have that
\begin{align}\label{kik1}
E(h_{n,r,l+1,A}(\mathcal{X}_{k+1})h_{n,r,l+1,A}(\mathcal{X}_{j}\cup (\mathcal{X}_{k+l+2-j}\backslash \mathcal{X}_{k+1})))\to c_r.
\end{align}
By combining equations  (\ref{sintag}), (\ref{kik2}) and (\ref{kik1}), and the definition of $\Phi_{j,A}(f_k.f_l)$, we conclude
\begin{align*}
&cov(f_{k,A}^\mathcal{P}(\rho),f_{l,A}^\mathcal{P}(\rho)) =\sum_{j=1}^{k}n^{k+l+2-j}r^{d(k+l+1-j)}\prod_{i=1}^l p_i^{\binom{k+1}{i+1}+ \binom{l+1}{i+1}- \binom{j}{i+1}}\Phi_{j,A}(f_k.f_l)
\end{align*}

\end{proof}
The following result is the the analogue  of Theorem~\ref{Candela13} for the soft random simplicial complexes $C(\mathcal{P}_n,r,\rho)$ and $C(n,r,\rho)$.
\begin{theorem}\label{Candela13x2}
Let $A\subset \mathbb{R}^d$ be an open set and $k,l\geq 1$. Consider the soft random simplicial complexes  $R(n,r,\rho)$ and $C(n,r,\rho)$ with the multiparameter $\rho=(p_1,p_2,\dots)$ and the Poisson versions $R(\mathcal{P}_n,r,\rho)$ and $C(\mathcal{P}_n,r,\rho)$.  Let  $f_{k,A}(\rho)$ and $f_{k,A}^\mathcal{P}(\rho)$ be the number of $k$-faces of $C(n,r,\rho)$ and $C(\mathcal{P}_n,r,\rho)$, respectively, whose left most point of the vertex set  is in $A$. 
For convenience, we define
$$\Theta_{j,A}(f_k,f_l)=\frac{\int_{A}f(x)^{k+l+2-j}dx}{j!(k+1-j)!(l+1-j)!}\times \int_{\mathbb{R}^d}...\int_{\mathbb{R}^d}\mathfrak{h}^j_{n,1,k+1,l+1}(0,x_2,...,x_{k+l+2-j})dx_2\dots dx_{k+l+j-2}.$$
Then, if $nr^d\to 0$, we have
\begin{align*}
cov(f_{k,A}^\mathcal{P}(\rho),f_{l,A}^\mathcal{P}(\rho))=\sum_{j=1}^{\min(k+1,l+1)}n^{k+l+2-j}r^{d(k+l+1-j)}\prod_{i=1}^l p_i^{\binom{k+1}{i+1}+ \binom{l+1}{i+1}- \binom{j}{i+1}}\Theta_{j,A}(f_k,f_l). 
\end{align*}
\end{theorem}

\begin{proof}
The proof for Theorem~\ref{Candela13x2} is the same as the proof of Theorem~\ref{Candela13}, replacing the functions $h$ with $\mathfrak{h}$, the constant $\Phi$ for $\Theta$ and  using Theorem~\ref{newcosito} instead of Proposition~\ref{Penrose-1}.
\end{proof}

\begin{corollary}\label{Candela14}
Consider the soft random simplicial complexes  $R(n,r,\rho)$ and $C(n,r,\rho)$ whose multiparameter $\rho=(p_1,p_2,\dots)$ and the Poisson versions $R(\mathcal{P}_n,r,\rho)$ and $C(\mathcal{P}_n,r,\rho)$. For $k\geq1$ let Let  $f_{k,A}(\rho)$ and $f_{k,A}^\mathcal{P}(\rho)$ be the number of $k$-faces of $C(n,r,\rho)$ and $C(\mathcal{P}_n,r,\rho)$, respectively, with the left most point of the vertex set  is in $A$, and suppose that $nr^d\to 0$. Then
$$var(f_{k,A}^\mathcal{P}(R(\mathcal{P}_n,r,\rho)))\sim \mu_{k+1,A}\prod_{i=1}^{k}p_i^{\binom{k+1}{i+1}}var(f_{k,A}^\mathcal{P}(R(\mathcal{P}_n,r)))$$
and
$$var(f_{k,A}^\mathcal{P}(C(\mathcal{P}_n,r,\rho)))\sim  \nu_{k+1,A}\prod_{i=1}^{k}p_i^{\binom{k+1}{i+1}}var(f_{k,A}^\mathcal{P}(C(\mathcal{P}_n,r)))$$
for $\mu_{k+1,A}$ and $\nu_{k+1,A}$ as in Definition~\ref{numudef}. 
\end{corollary}
\begin{proof} By Theorem~\ref{Candela13}, we have
\begin{align}\label{ala1}
var(f_{k,A}^\mathcal{P}(R(\mathcal{P}_n,r,\rho)))&=cov(f_{k,A}^\mathcal{P}(R(\mathcal{P}_n,r,\rho)),f_{k,A}^\mathcal{P}(R(\mathcal{P}_n,r,\rho)))\\
\notag&=\sum_{j=1}^{k+1}n^{2k+2-j}r^{d(2k+1-j)}\prod_{i=1}^k p_i^{2\binom{k+1}{i+1}- \binom{j}{i+1}}\Phi_{j,A}(f_k,f_k)
\end{align}
and by Theorem~\ref{Candela13x2},
\begin{align}\label{ala2}
var(f_{k,A}^\mathcal{P}(C(\mathcal{P}_n,r,\rho)))&=cov(f_{k,A}^\mathcal{P}(C(\mathcal{P}_n,r,\rho)),f_{k,A}^\mathcal{P}(C(\mathcal{P}_n,r,\rho)))\\
\notag&=\sum_{j=1}^{k+1}n^{2k+2-j}r^{d(2k+1-j)}\prod\limits_{i=1}^k p_i^{2\binom{k+1}{i+1}- \binom{j}{i+1}}\Theta_{j,A}(f_k,f_k).
\end{align}
Comparing the $(j+1)$-st term with the
$j$-th term in the expressions above, we obtain
$$\frac{n^{2k+2-j}r^{d(2k+1-j)}\prod\limits_{i=1}^{k} p_i^{2\binom{k+1}{i+1}- \binom{j}{i+1}}}{n^{2k+2-j-1}r^{d(2k+1-j-1)}\prod\limits_{i=1}^{k} p_i^{2\binom{k+1}{i+1}-\binom{j+1}{i+1}}}=\frac{n^{2k+2-j}r^{d(2k+1-j)}\prod\limits_{i=1}^{j-1} p_i^{- \binom{j}{i+1}}}{n^{2k+2-j-1}r^{d(2k+1-j-1)}\prod\limits_{i=1}^{j} p_i^{-\binom{j+1}{i+1}}}= nr^d \prod\limits_{i=1}^{j} p_i^{ \binom{j}{i}}\to 0.$$
If we divide Equations (\ref{ala1}) and (\ref{ala2}) by the term $n^{k+1}r^{dk}\prod\limits_{i=1}^{k+1} p_i^{\binom{k+1}{i+1}}$, we get
\begin{align*}
\frac{var(f_{k,A}^\mathcal{P}(R(\mathcal{P}_n,r,\rho)))}{n^{k+1}r^{dk}\prod\limits_{i=1}^{k+1} p_i^{\binom{k+1}{i+1}}}&=\frac{cov(f_{k,A}^\mathcal{P}(R(\mathcal{P}_n,r,\rho)),f_{k,A}^\mathcal{P}(C(\mathcal{P}_n,r,\rho)))}{n^{k+1}r^{dk}\prod\limits_{i=1}^{k+1} p_i^{\binom{k+1}{i+1}}}\sim\Phi_{k+1,A}(f_k,f_k).
\end{align*}
and
\begin{align*}
\frac{var(f_{k,A}^\mathcal{P}(C(\mathcal{P}_n,r,\rho)))}{n^{k+1}r^{dk}\prod\limits_{i=1}^{k+1} p_i^{\binom{k+1}{i+1}}}&=\frac{cov(f_{k,A}^\mathcal{P}(C(\mathcal{P}_n,r,\rho)),f_{k,A}^\mathcal{P}(C(\mathcal{P}_n,r,\rho)))}{n^{k+1}r^{dk}\prod\limits_{i=1}^{k+1} p_i^{\binom{k+1}{i+1}}}\sim\Theta_{k+1,A}(f_k,f_k).
\end{align*}
Also, notice that for $\mu_{k+1,A}$ and $\nu_{k+1,A}$ as in Definition~\ref{numudef}, $\Theta_{k+1,A}(f_k,f_k)=\nu_{k+1,A}$ and $\Phi_{k+1,A}(f_k,f_k)=\mu_{k+1,A}$. Thus
\begin{align*}
var(f_{k,A}^\mathcal{P}(R(\mathcal{P}_n,r,\rho)))\sim \mu_{k+1,A} n^{k+1}r^{dk}\prod_{i=1}^k p_i^{\binom{k+1}{i+1}}\\
var(f_{k,A}^\mathcal{P}(C(\mathcal{P}_n,r,\rho)))\sim \nu_{k+1,A}n^{k+1}r^{dk}\prod_{i=1}^k p_i^{\binom{k+1}{i+1}}. 
\end{align*}
Furthermore, notice that we can recover $R(\mathcal{P}_n,r)$ and $C(\mathcal{P}_n,r)$ from the soft random simplicial complexes $R(\mathcal{P}_n,r,\rho))$ and $R(\mathcal{P}_n,r,\rho))$  by setting each $p_i=1$, thus  by Theorems~\ref{Candela13} and \ref{Candela13x2}
$$var(f_{k,A}^\mathcal{P}(R(\mathcal{P}_n,r)))\sim \mu_{k+1,A}n^{k+1}r^{dk}$$
and
$$var(f_{k,A}^\mathcal{P}(C(\mathcal{P}_n,r)))\sim \nu_{k+1,A}n^{k+1}r^{dk}.$$
Therefore, we have that:
$$var(f_{k,A}^\mathcal{P}(R(\mathcal{P}_n,r,\rho)))\sim\mu_{k+1,A}\prod_{i=1}^{k}p_i^{\binom{k+1}{i+1}} n^{k+1}r^{dk} \sim\prod_{i=1}^{k}p_i^{\binom{k+1}{i+1}}var(f_{k,A}^\mathcal{P}(R(\mathcal{P}_n,r)))$$
and
$$var(f_{k,A}^\mathcal{P}(C(\mathcal{P}_n,r,\rho)))\sim\nu_{k+1,A}\prod_{i=1}^{k}p_i^{\binom{k+1}{i+1}} n^{k+1}r^{dk} \sim\prod_{i=1}^{k}p_i^{\binom{k+1}{i+1}}var(f_{k,A}^\mathcal{P}(C(\mathcal{P}_n,r))).$$
\end{proof}

\begin{theorem}\label{Candela10}
Consider the soft random simplicial complexes $R(n,r,\rho)$ and $C(n,r,\rho)$ with multiparameter $\rho=(p_1,p_2,\dots)$ and the Poisson versions $R(\mathcal{P}_n,r,\rho)$ and $C(\mathcal{P}_n,r,\rho)$. Suppose that $nr^d\to 0$ and  $k,l\geq1$. Then

$$\frac{var(f_{k,A}^\mathcal{P}(\rho))}{var(f_{0,A}^\mathcal{P}(\rho))}\to 0,~~~~\frac{cov(f_{k,A}^\mathcal{P}(\rho),f_l^\mathcal{P}(\rho))}{var(f_{0,A}^\mathcal{P}(\rho))}\to 0.$$
\end{theorem}
\begin{proof} 
Notice that in both cases,  number of vertices of the complex $f_{0,A}^{\mathcal{P}}(\rho)$ does not depend on the multiparameter $\rho=(p_1,p_2,...,p_{n-1}, \dots)$. Thus, $f_{0,A}^{\mathcal{P}}(\rho)$ is the number of vertices in the geometric graph $G(\mathcal{P}_n,r)$ that lie in $A$. Therefore,  $f_{0,A}^{\mathcal{P}}(\rho)$ is a Poisson random variable with rate $\Lambda$ for
$$\Lambda= \int_{A}nf(x)dx\sim \mu_{0,A} n$$
with
\begin{equation}\label{mu0}
\mu_{0,A}= \nu_{0,A}=\int_{A}f(x)dx.
\end{equation}
It follows that 
\begin{equation}
\label{eq:Variance}
var(f_{0,A}^\mathcal{P}(\rho))=\mu_{0,A}n.
\end{equation} 
By Corollary~\ref{Candela14}, for both simplicial complexes it is true that, for a certain constant $c$,
$$\frac{var(f_k^\mathcal{P}(\rho))}{ var(f_0^\mathcal{P}(\rho))}\sim \frac{c \prod\limits_{i=1}^{k}p_i^{\binom{k+1}{i+1}}n^{k+1}r^{dk}}{\mu_{0,A}n}\leq \frac{c}{\mu_{0,A}}(nr^{d})^k \to 0.$$
On the other hand, by Theorem~\ref{Candela13} and by Theorem~\ref{Candela13x2}, for a certain constant $c_j$ depending on $j$,
$$\frac{cov(f_k^\mathcal{P}(\rho),f_l^\mathcal{P}(\rho))}{var(f_0^\mathcal{P}(\rho))}\sim \frac{\sum\limits_{j=1}^{k}c_jn^{k+l+2-j}r^{d(k+l+1-j)}\prod\limits_{i=1}^l p_i^{\binom{k+1}{i+1}+ \binom{l+1}{i+1}- \binom{j}{i+1}}}{\mu_{0,A}n} $$
$$\leq  \frac{\max_{j}c_j}{\mu_{0,A}}\sum_{j=1}^k {(nr^d)}^{k+l+1-j}\to 0.$$
\end{proof}

\begin{theorem}\label{varf0}
Consider the soft random simplicial complexes  $R(n,r,\rho)$ and $C(n,r,\rho)$ with the multiparameter $\rho=(p_1,p_2,...)$ and the Poisson versions $R(\mathcal{P}_n,r,\rho)$ and $C(\mathcal{P}_n,r,\rho)$. Suppose that $nr^d\to 0$ and  $l\geq 1$. Then
$$\frac{cov(f_{l,A}^\mathcal{P}(\rho),f_0^\mathcal{P}(\rho))}{var(f_{0,A}^\mathcal{P}(\rho))}\to 0.$$
\end{theorem}

\begin{proof}
We will prove the theorem for the soft random simplicial complex $R(\mathcal{P}_n,r,\rho)$ (the proof for $C(\mathcal{P}_n,r,\rho)$ is analogous, replacing all instances of $h$ with  $\mathfrak{h}$ and  using Theorem~\ref{newcosito} instead of Proposition~\ref{Penrose-1}.
We first observe that
\begin{align*}
cov(f_{l,A}^\mathcal{P}(\rho),f_{0,A}^\mathcal{P}(\rho))&=E(f_{l,A}^\mathcal{P}(\rho)f_{0,A}^\mathcal{P}(\rho))-E(f_{l,A}^\mathcal{P}(\rho))E(f_{0,A}^\mathcal{P}(\rho))\\
&=E\left(\sum_{j=0}^1\sum_{Y'\subset \mathcal{P}_n}\sum_{y\in \mathcal{P}_n}h_{\mathcal{P}_n,r,l+1,A,\rho}(Y')\cdot 1_{y \in A}(y)1_{\{|\{y\}\cap Y'|=j\}}\right)\\
&\qquad-E(f_{l,A}^\mathcal{P}(\rho))E(f_0^\mathcal{P}(\rho)).
\end{align*}
 By Theorem~\ref{Penrose-p}, the term $j=0$ on the left-hand side of the expression above corresponds to $E(f_{0,A}^\mathcal{P}(\rho))E(f_{l,A}^\mathcal{P}(\rho))$, which cancels the term on the right. We therefore have,
\begin{align*}
cov(f_{l,A}^\mathcal{P}(\rho),f_{0,A}^\mathcal{P}(\rho))&=E\left(\sum_{Y'\subset \mathcal{P}_n}\sum_{y\in \mathcal{P}_n}h_{\mathcal{P}_n,r,l+1,A,\rho}(Y')1_{y \in A}(y)1_{y \in Y'}(y)\right)\\
&=E\left(\sum_{X\subset \mathcal{P}_n}H(X)\right),
\end{align*}
where 
$$H(X)=1_{|X|=l+1}\sum_{y\in X}h_{\mathcal{P}_n,r,l+1,A,\rho}(X)1_{y \in A}(y)1_{y \in X}(y).$$
By Theorem~\ref{Penrose-5} (Palm theory),
\begin{align*}
cov(f_{l,A}^\mathcal{P}(\rho),f_0^\mathcal{P}(\rho))&=\frac{n^{l+1}}{(l+1)!}E(H(\mathcal{X}_{l+1}))\\ 
&= \frac{n^{l+1}}{(l+1)!}(l+1)E(h_{n,r,l+1,A,\rho}(\mathcal{X}_l\cup\{y\})\cdot1_{y \in A}(y) ),
\end{align*}
where $y \in \mathcal{X}_{l+1}$. By Theorem~\ref{prodhx1}, this expression becomes
\begin{align*}
&\qquad\qquad\qquad\qquad\qquad= \frac{n^{l+1}}{l!}\prod_{i=1}^{l}p_i^{\binom{l+1}{i+1}}E(h_{n,r,l+1,A}(\mathcal{X}_l\cup\{y\}) \cdot1_{y \in A}(y) ).    
\end{align*}
Thus, by Proposition~\ref{Penrose-1} and Theorem~\ref{equalh},
$$cov(f_{l,A}^\mathcal{P}(\rho),f_{0,A}^\mathcal{P}(\rho))\sim   c n^{l+1}r^{dl}\prod_{i=1}^{l}p_i^{\binom{l+1}{i+1}}.$$
Furthermore, by Equation \ref{eq:Variance}, $var(f_{0,A}^\mathcal{P}(\rho))=\mu_{0,A}n$, and therefore
$$ \frac{cov(f_{l,A}^\mathcal{P}(\rho),f_0^\mathcal{P}(\rho))}{var(f_{0,A}^\mathcal{P}(\rho))}\sim \frac{c n^{l}r^{dl}\prod\limits_{i=1}^{l}p_i^{\binom{l+1}{i+1}}}{\mu_{0,A}}\to 0.$$

\end{proof}

\begin{lemma}\label{Candela11}
Consider the soft random simplicial complexes $R(n,r,\rho)$ and $C(n,r,\rho)$ with multiparameter $\rho=(p_1,p_2,...,p_{n-1}, \dots)$ and the Poisson versions $R(\mathcal{P}_n,r,\rho)$ and $C(\mathcal{P}_n,r,\rho)$. Let $\Sigma_1,\Sigma_2$ simplicial complexes such that $\Sigma_1\in\{R(\mathcal{P}_n,r,\rho), C(\mathcal{P}_n,r,\rho)\}$ and $\Sigma_2\in\{R(n,r,\rho), C(n,r,\rho)\}$.
For the simplicial complexes $\Sigma_1$ and $\Sigma_2$, define  $\chi_{A}^\mathcal{P}(\rho,\Sigma_1)$  and $\chi_{A}(\rho,\Sigma_2)$ as:
$$\chi_{A}^\mathcal{P}(\rho,\Sigma_1)=\sum_{j=0}^{n-1}(-1)^j f_{j,A}^\mathcal{P}(\rho, \Sigma_1).$$
$$\chi_{A}(\rho, \Sigma_2)=\sum_{j=0}^{n-1}(-1)^j f_{j,A}(\rho,\Sigma_2).$$
Then, under the hypothesis of Theorem~\ref{Candela16}, we have 
$$\chi_{A}^\mathcal{P}(\rho,\Sigma_1)=\sum_{j=0}^{l}(-1)^j f_{j,A}^\mathcal{P}(\rho, \Sigma_1)~~~\text{a.a.s}$$
and
$$\chi_{A}(\rho,\Sigma_2)=\sum_{j=0}^{l}(-1)^j f_{j,A}(\rho, \Sigma_2)~~~a.a.s.$$
Furthermore, 
$$var(\chi_{A}^\mathcal{P}(\rho, \Sigma_1))\sim var(f_0^\mathcal{P}(\rho))~~~\text{a.a.s}$$
\end{lemma}

\begin{proof}
We will prove the theorem for the soft random simplicial complexes $R(n,r,\rho)$and $R(\mathcal{P}_n,r,\rho)$ (the proof for $C(n,r,\rho)$ and $C(\mathcal{P}_n,r,\rho)$ is analogous, replacing all instances of $h$ for $\mathfrak{h}$).\\
Let us start first with $R(n,r,\rho)$. The expected number of $(l+1)$-faces of  $R(n,r,\rho)$ with left most point in $A$ satisfies
\begin{align*}
E(f_{l+1,A}(\rho))&=E\left(\sum_{X\subset \mathcal{X}_n } h_{n,r,l+2,A,\rho}(X)\right)\\
&=\binom{n}{l+2}E(h_{n,r,l+2,A,\rho}(\mathcal{X}_{l+2}))    
\end{align*}
by Theorem~\ref{hvsh},
\begin{align*}
&\label{eq:Soft to hard}\qquad\qquad\qquad\qquad=\binom{n}{l+2}\prod_{i=1}^{l+1}p_i^{\binom{l+2}{i+1}}E(h_{n,r,l+2,A}(\mathcal{X}_{l+2}))\\
&\qquad\qquad\qquad\qquad=\prod_{i=1}^{l+1}p_i^{\binom{l+2}{i+1}} E(f_{l+1,A}),
\end{align*}
where $E(f_{l+1,A})$ is the expected number of $(l+1)$-faces in $R(n,r)$ with left most point of the $l+1$-face in $A$. On the other hand, for $R(\mathcal{P}_n,r,\rho)$ we have:
$$E(f_{l+1,A}^\mathcal{P}(\rho))=E\left(\sum_{Y\subset \mathcal{P}_n}h_{\mathcal{P}_n,r,l+2,A,\rho}(Y)\right)$$
using Theorem~\ref{Penrose-5} (Palm theory), this results in
\begin{align*}
&\qquad\qquad\qquad=\frac{n^{l+2}}{(l+2)!}E(h_{n,r,l+2,A,\rho}(\mathcal{X}_{l+2})) .
\end{align*}
Furthermore, by Theorem~\ref{hvsh}, this equals
\begin{align*}
&\qquad\qquad\qquad\qquad=\frac{n^{l+2}}{(l+2)!}\prod_{i=1}^{l+1}p_i^{\binom{l+2}{i+1}}E(h_{n,r,l+2,A}(\mathcal{X}_{l+2}))\\
&\qquad\qquad\qquad\qquad\sim \prod_{i=1}^{l+1}p_i^{\binom{l+2}{i+1}} \binom{n}{l+2} E(h_{n,r,l+2,A}(\mathcal{X}_{l+2}))\\
&\qquad\qquad\qquad\qquad\sim\prod_{i=1}^{l+1}p_i^{\binom{l+2}{i+1}} E(f_{l+1,A}).   
\end{align*}
Therefore, in both cases: 
$$E(f_{l+1,A}^{\mathcal{P}}(\rho))\sim E(f_{l+1,A}(\rho))\sim \prod_{i=1}^{l+1}p_i^{\binom{l+2}{i+1}} E(f_{l+1,A}).$$
On the other hand, by Proposition~\ref{Penrose-1} and Theorem~\ref{equalh}, $E(f_{l+1, A})\sim \mu_{l+2,A}n^{l+2}r^{d(l+1)}$. Thus, by Markov's inequality, Lemma \ref{Frieze2}, we have 
$$P(f_{l+1,A}^\mathcal{P}(\rho)>0)\leq E(f_{l+1,A}^\mathcal{P}(\rho))\sim \prod_{i=1}^{l+1}p_i^{\binom{l+2}{i+1}}\mu_{l+2,A}n^{l+2}r^{d(l+1)}\to 0$$
and
$$P(f_{l+1,A}(\rho)>0)\leq E(f_{l+1,A}(\rho))\sim \prod_{i=1}^{l+1}p_i^{\binom{l+2}{i+1}}\mu_{l+2,A}n^{l+2}r^{d(l+1)}\to 0$$
which tend to zero since by hypothesis $\prod\limits_{i=1}^{l+1}p_i^{\binom{l+2}{i+1}}n^{l+2}r^{d(l+1)}\to 0$ and therefore $f_{l+1,A}(\rho)=0$ and $f_{l+1,A}^\mathcal{P}(\rho)=0$, a.a.s. Furthermore, for every $k\geq{l+1}$, we have
$$E(f_{k,A}(\rho))\sim \prod_{i=1}^{k}p_i^{\binom{k+1}{i+1}} \mu_{k+1,A}n^{k+1}r^{dk}=\prod_{i=1}^{k}p_i^{\binom{k+1}{i+1}}n^{l+2}\mu_{k+1,A}r^{d(l+1)} (nr^{d})^{k-l-1}\to 0$$
and 
$$E(f_{k,A}^\mathcal{P}(\rho))\sim \prod_{i=1}^{k}p_i^{\binom{k+1}{i+1}} \mu_{k+1,A}n^{k+1}r^{dk}=\prod_{i=1}^{k}p_i^{\binom{k+1}{i+1}}n^{l+2}\mu_{k+1,A}r^{d(l+1)} (nr^{d})^{k-l-1}\to 0$$
which tend to zero, since, by hypothesis, $\prod\limits_{i=1}^{l+1}p_i^{\binom{l+2}{i+1}}n^{l+2}r^{d(l+1)}\to 0$ and $nr^d\to 0$. Therefore, by Markov's inequality, Lemma \ref{Frieze2}, for every $k\geq{l+1}$,  $f_{k,A}^\mathcal{P}(\rho)=f_{k,A}(\rho)=0$ a.a.s.

Finally, given the above, it follows that a.a.s
$$\chi_{A}^\mathcal{P}(\rho)=\sum_{j=0}^{l}(-1)^j f_{j,A}^\mathcal{P}(\rho) ~~~~~~~\text{      and      }~~~~~~~\chi_{A}(\rho)=\sum_{j=0}^{l}(-1)^j f_{j,A}(\rho),~~~\text{a.a.s}$$
and therefore, a.a.s by properties of the variance 
$$var(\chi^\mathcal{P}_A(\rho))=\sum_{i=0}^l var(f_i^\mathcal{P}(\rho))+\sum_{i,j}(-1)^{i+j}cov(f_{i,A}^\mathcal{P}(\rho),f_{j,A}^\mathcal{P}(\rho)).$$
thus:
$$\frac{var(\chi_A^\mathcal{P}(\rho))}{var(f_{0,A}^\mathcal{P}(\rho))}=\sum_{i=0}^l \frac{var(f_{i,A}^\mathcal{P}(\rho))}{var(f_{0,A}^\mathcal{P}(\rho))}+\sum_{i,j}(-1)^{i+j}\frac{cov(f_i^\mathcal{P}(\rho),f_{j,A}^\mathcal{P}(\rho))}{var(f_{0,A}^\mathcal{P}(\rho))}$$
and therefore, by Theorems~\ref{Candela10} and \ref{varf0}:
$$\frac{var(\chi_A^\mathcal{P}(\rho))}{var(f_{0,A}^\mathcal{P}(\rho))}\sim 1.$$
\end{proof}

\begin{lemma}\label{Hyj}
Let $X\geq 0$ be a random variable. For $q=3,4$, there exists a constant $c_q$ such that
$$E|X- E(X)|^q\leq c_q E|X|^q$$
\end{lemma}
\begin{proof}
Let us start with $q=3$. Notice that
$$|X- E(X)|^3\leq (X+E(X))^3= X^3 +3XE(X)^2+2X^2E(X)+E(X)^3$$
then
\begin{equation}\label{009}
E|X- E(X)|^3\leq E(X^3)+2E(X^2)E(X)+4E(X)^3. 
\end{equation}
By H\H{o}lder's inequality, \ref{Holder}, we have
$$E(X^2)\leq E(X^3)^{\frac{2}{3}} $$
and
$$E(X)\leq E(X^3)^{\frac{1}{3}}.$$
On the other hand, since $X\geq 0$, $X^3$ is a convex function and therefore by Jensen's inequality, Lemma~\ref{Jensen}, we have
$$E(X)^3\leq E(X^3).$$
Thus, combining equation~(\ref{009}) with Jensen's and  H\H{o}lder,s inequality, we have
$$E|X- E(X)|^3\leq 7E(X^3).$$
On the other hand, for $q=4$, notice that
$$|X- E(X)|^4\leq (X+E(X))^4= X^4 +4XE(X)^3+6X^2E(X)^2+4X^3E(X)+E(X)^4$$
then
\begin{equation}\label{0091}
E|X- E(X)|^4\leq E(X^4)+6E(X^2)E(X)^2+4E(X^3)E(X)+5E(X)^4
\end{equation}
Now, by Jensen's inequality, Lemma~\ref{Jensen}, we have 
$$E(X)^2\leq E(X^2)$$
and
$$E(X)^4\leq E(X^4).$$
on the other hand, by H\H{o}lder's inequality, Lemma~\ref{Holder}, we have
$$E(X^2)\leq E(X^4)^{\frac{2}{4}} $$
also,
$$E(X^3)\leq E(X^4)^{\frac{3}{4}} $$
and
$$E(X)\leq E(X^4)^{\frac{1}{4}}.$$
Thus, combining equation Equation~(\ref{0091}) with Jensen's and  H\H{o}lder,s inequality, we have
$$E|X- E(X)|^4\leq 16 E(X^4).$$

\end{proof}

\begin{lemma}\label{pq}

Let $\{Q_{m,n}\}_{m\in \mathbb{N}}$ be a partition of $\mathbb{R}^d$ into cubes of side length $r$ and let $A\subset \mathbb{R^d}$ be open and bounded.
Let $\xi_{m,\rho}$ counts the number of $k$-faces of $R(\mathcal{P}_n,r,\rho)$ with left most point of the $k$-face in $A\cap Q_{m,n}$ and   $\xi_{m}$ is number of $k$-faces of $R(\mathcal{P}_n,r)$ with left most point in $A\cap Q_{m,n}$. Then, for $q=3,4$ , we have
\begin{equation}
E|\xi_{m,\rho}|^q\leq \prod\limits_{i=1}^k p_i^{\binom{k+1}{i+1}} E|\xi_{m}|^q    
\end{equation}
\end{lemma}
\begin{proof}
we first start with $q=3$.
\begin{align}\label{plpl}
 E|\xi_{i,\rho}|^3&=E\left(\sum_{Y_1,Y_2,Y_3\subset \mathcal{P}_n} h_{\mathcal{P}_n,r,k+1,B,\rho}(Y_1)h_{\mathcal{P}_n,r,k+1,(A\cap Q_{m,n}),\rho}(Y_2)h_{\mathcal{P}_n,r,k+1,(A\cap Q_{m,n}),\rho}(Y_3)\right).
\end{align}
Consider now the variables $s=|Y_1\cap Y_2\cap Y_3|$, $j_1=|Y_1\cap Y_2\cap Y_3^c|$, $j_2=|Y_1^c\cap Y_2\cap Y_3|$ and $j_3=|Y_1\cap Y_2^c\cap Y_3|$. Notice that, by definition, $s,j_1,j_2,j_3$ satisfy that
$$j_1+j_2+s\leq k+1$$
$$j_3+s\leq k+1$$
and therefore,
\begin{align}\label{sj}
3k+3-j_1-j_2-j_3-2s\geq k+1    
\end{align}
On the other hand, define  the functions $Q(Y_1,Y_2,Y_3)$ and $S(j_1,j_2,j_3,s)$ as 
$$Q(Y_1,Y_2,Y_3)=h_{\mathcal{P}_n,r,k+1,(A\cap Q_{m,n}),\rho}(Y_1)h_{\mathcal{P}_n,r,k+1,(A\cap Q_{m,n}),\rho}(Y_2)h_{\mathcal{P}_n,r,k+1,(A\cap Q_{m,n}),\rho}(Y_3)$$
and
$$S(j_1,j_2,j_3,s)=1_{\{|Y_1\cap Y_2\cap Y_3|=s\}\}}1_{\{|Y_1\cap Y_2\cap Y_3^c|=j_1\}\}}1_{\{|Y_1^c\cap Y_2\cap Y_3|=j_2\}\}}1_{\{|Y_1\cap Y_2^c\cap Y_3|=j_3\}\}}$$
Then, from Equation (\ref{plpl}), we have that 
\begin{align*}
E|\xi_{m,\rho}|^3&=E\left(\sum_{s=0}^{k+1}\sum_{j_1,j_2,j_3=0}^{k+1-s}\sum_{Y_1,Y_2,Y_3\subset \mathcal{P}_n} Q(Y_1,Y_2,Y_3)S(j_1,j_2,j_3,s)\right)\\
&= E\left(\sum_{s=0}^{k+1}\sum_{j_1,j_2,j_3=0}^{k+1-s}\sum_{X\subset \mathcal{P}_n} H_{j_1,j_2,j_3,s}(X)\right)
\end{align*}
for
$$H_{j_1,j_2,j_3,s}(X)=1_{\{|X|=3k+3-j_1-j_2-j_3-s\}}\sum_{Y_1,Y_2,Y_3\subset X}Q(Y_1,Y_2,Y_3)S(j_1,j_2,j_3,s).$$
Using Theorem~\ref{Penrose-5} (Palm theory),  applied to $H_{j_1,j_2,j_3,s}(X)$

\begin{align*}
E\left(\sum_{X\subset \mathcal{P}_n} H_{j_1,j_2,j_3,s}(X)\right)&= \frac{n^{3k+3-j_1-j_2-j_3-2s}}{(3k+3-j_1-j_2-j_3-2s)!}\\
&\times E(h_{n,r,k+1,(A\cap Q_{m,n}),\rho}(\mathcal{X}_{k+1})h_{n,r,k+1,(A\cap Q_{m,n}),\rho}(\mathcal{X}'_{k+1})h_{n,r,k+1,(A\cap Q_{m,n}),\rho}(\mathcal{X}''_{k+1})   
\end{align*}
where the sets $\mathcal{X}_{k+1},\mathcal{X}'_{k+1},\mathcal{X}''_{k+1}\subset \mathcal{X}'_{n}$ have cardinality $k+1$ and satisfy that
$$s=|\mathcal{X}_{k+1}\cap\mathcal{X}'_{k+1}\cap \mathcal{X}''_{k+1}|$$ 
$$j_1=|\mathcal{X}_{k+1}\cap \mathcal{X}'_{k+1}\cap (\mathcal{X}''_{k+1})^c|$$
$$j_2=|(\mathcal{X}_{k+1})^c\cap \mathcal{X}'_{k+1}\cap \mathcal{X}''_{k+1}|$$ 
and 
$$j_3=|\mathcal{X}_{k+1}\cap (\mathcal{X}'_{k+1})^c\cap \mathcal{X}''_{k+1}|.$$
Therefore, by Theorem~\ref{prodh's}, the expression above is equal to 
\begin{align*}
&= E(\sum_{j=0}^{k+1} \frac{n^{3k+3-j_1-j_2-j_3-2s}}{(3k+3-j_1-j_2-j_3-2s)!}\prod\limits_{i=1}^k p_i^{3\binom{k+1}{i+1}-\binom{j_1}{i+1}-\binom{j_2}{i+1}-\binom{j_3}{i+1}-\binom{2s}{i+1}}\\
&\times E(h_{n,r,k+1,(A\cap Q_{m,n})}(\mathcal{X}_{k+1})h_{n,r,k+1,(A\cap Q_{m,n})}(\mathcal{X}'_{k+1})h_{n,r,k+1,(A\cap Q_{m,n})}(\mathcal{X}''_{k+1})  
\end{align*}
and by Equation~(\ref{sj}), we have
\begin{align*}
&\leq \prod\limits_{i=1}^k p_i^{\binom{k+1}{i+1}} E(\sum_{j=0}^{k+1} \frac{n^{3k+3-j_1-j_2-j_3-2s}}{(3k+3-j_1-j_2-j_3-2s)!}\\
&\times E(h_{n,r,k+1,(A\cap Q_{m,n})}(\mathcal{X}_{k+1})h_{n,r,k+1,(A\cap Q_{m,n})}(\mathcal{X}'_{k+1})h_{n,r,k+1,(A\cap Q_{m,n})}(\mathcal{X}''_{k+1})  
\end{align*}
now, by returning over the same argument
$$E|\xi_{m,\rho}|^3\leq \prod\limits_{i=1}^k p_i^{\binom{k+1}{i+1}} E|\xi_{m}|^3.$$ 
The proof for $q=4$ is analogous.

\end{proof}

The following Lemma and its proof were taken from the proof of Theorem 3.15 in \cite{MR3079211}.
\begin{lemma}\label{z}
Let $\{Q_{m,n}\}_{m\in \mathbb{N}}$ be a partition of $\mathbb{R}^d$ into cubes of side length $r$ and let $A\subset \mathbb{R^d}$ be open and bounded.
Let $\xi_{m}$ be the number of $k$-faces of $R(\mathcal{P}_n,r)$ with left most point in a bounded set $A\cap Q_{m,n}$ and let  $z_m$ be the number of points of $\mathcal{P}_n$ within a  distance $2r$ of the cube $Q_{m,n}$, then for  $(z_m)_{k}=z_i({z_m}-1)...(z_m-j)$, we have
$$|\xi_{m}|^q\leq (z_m)_{k}=z_m({z_m}-1)...(z_m-j).$$
Therefore, for a constant $c_k$
$$ E|\xi_{m}|^q\leq E((z_m)_k^q)\leq c_k(nr^d)^{k+1}.$$
\end{lemma}

\begin{proof}
observe that if  $z_i$ is the number of points of $\mathcal{P}_n$ within a  distance $2r$ of the cube $Q_{m,n}$ then $z_m$ is distributed as a Poisson random variable with mean $\lambda$ for 
$$\lambda= n\int\limits_{B(Q_{m,n},2r)}f(x)dx$$
therefore, there exists a constant $c$ such that
$$\lambda= n\int\limits_{B(Q_{m,n},2r)}f(x)dx\leq cnr^d.$$
Now, writing $(z_m)_{k}=z_m({z_m}-1)...(z_m-j)$, it follows that
$$|\xi_{m}|^q\leq (z_m)_{k}=z_m({z_m}-1)...(z_m-j)$$
and therefore,
\begin{equation}\label{bbb}
E|\xi_{m}|^q\leq E((z_m)_{k})\leq \sum_{l=k+1}^{\infty}(l)_k^q\frac{e^{-cnr^d} (cnr^d)^l}{l!}.    
\end{equation}
Comparing the terms $l-th$ and $(l+1)$-th in the sum above, we have
$$\frac{\frac{(l+1)_k^qe^{-cnr^d} (cnr^d)^{l+1}}{(l+1)!}}{\frac{(l)_k^qe^{-cnr^d} (cnr^d)^l}{l!}}=\frac{l!(l+1)_k^q(cnr^d)^{l+1}}{(l+1)!(l)_k^q(cnr^d)^{l}}=\frac{(l+1)^{q-1}cnr^d}{(l-k)^q}\to 0,$$
and therefore dividing Equation~(\ref{bbb}) by $\frac{(k+1)_k^qe^{-cnr^d} (cnr^d)^{k+1}}{(k+1)!}$, we have 
$$\frac{E|\xi_{m}|^q}{\frac{(k+1)_k^qe^{-cnr^d} (cnr^d)^{k+1}}{(k+1)!}}\leq \sum_{l=k+1}^{\infty}\frac{\frac{(l)_k^qe^{-cnr^d} (cnr^d)^l}{l!}}{\frac{(k+1)_k^q e^{-cnr^d} (cnr^d)^{k+1}}{(k+1)!}}\to 1,$$
it follows that, asymptotically, 
$$E|\xi_{m}|^q\leq \frac{(k+1)_k^qe^{-cnr^d} (cnr^d)^{k+1}}{(k+1)!} $$
Thus, asymptotically, there exists a constant $c_k$ such that
$$E|\xi_{m}|^q\leq E((z_m)_{k})\leq c_kn^{k+1}r^{d(k+1)}. $$
\end{proof}

\section{C.L.T for the count $f_k$ of $R(\mathcal{P}_n,r,\rho)$ and $C(\mathcal{P}_n,r,\rho)$}\label{section4}

In order to prove the Theorems~\ref{Candela15} and \ref{Candela16}, we will first consider the Poisson version of the problem. That is, we first prove each version of the central limit theorem for  the soft random simplicial complexes $R(\mathcal{P}_n,r,\rho)$ and $C(\mathcal{P}_n,r,\rho)$, and we then  prove  the  original theorem using a depoissonization argument. The principal reason for approaching the problem in this way  is to use  the strong spatial independence properties of the Poisson process $\mathcal{P}_n$ together with a normal approximation theorem (Theorem~\ref{Penrose-4}) to  facilitate the argument in the  Poisson version of the problem, from which we may then  recover the result for the original model  by a depoissonization argument.  This technique is very common when dealing with central limit theorems for  random graphs or random simplicial complexes built on a  binomial process $\mathcal{X}_n$. The reader may also consult \cite{MR1986198} and \cite{MR3079211} to see this strategy used to prove a variety of other central limit theorems.

We now quote from \cite{MR1986198}  the normal approximation theorem that we will use to prove the central limit theorems for the Poisson problem. We begin with the following definition.

\begin{definition}\label{dependency}
Suppose $(I,E)$ is a graph with finite or countable vertex set $I$. For $i,j\in I$ write $i\sim j$ iff $\{i,j\}\in E$. For for $i\in I$, let $N_i$ denote the adjacency neighborhood of $i$, that is, the set $\{i\}\cup\{j\in I, \}$.  We say that the graph $(I,\sim)$ is a dependency graph for a collection of random variables $\{\zeta_i\}_{i\in I}$ if for any two disjoints subsets $I_1,I_2$ of $I$ such that there are no edges connecting $I_1$ with $I_2$, the collection of random variables $\{\zeta_i\}_{i\in I_1}$ is independent of $\{\zeta_i\}_{i\in I_2}.$
\end{definition}

\begin{theorem}[Theorem 2.4 \cite{MR1986198}]\label{Penrose-4} Suppose $\{\zeta_i\}_{i\in I}$ is a finite collection of random variables with $E(\xi_i)=0$
 for each $i$, and whose dependency graph $(I,\sim)$ has maximum degree $D-1$. Set $W=\sum_{i\in I}\zeta_i$ and suppose $E(W^2)=1$. Let $Z$ be a standard normal random variable. Then for all $t\in \mathbb{R}$,
$$|P(W\leq t)-P(Z\leq t)|\leq \frac{2}{\sqrt[4]{2\pi}}\sqrt{D^2\sum_{i\in I}E|\zeta_i|^3}+ 6 \sqrt{D^3\sum_{i\in I}E|\zeta_i|^4} $$ 
\end{theorem}

We present now the Poisson version of Theorem~\ref{Candela15} for bounded sets $A\subset \mathbb{R}^d$,  the central limit theorem for the count of $k$-faces $f_{k,A}^{\mathbb{P}}$ in the simplicial complexes $R(\mathcal{P}_n,r,\rho)$ and $C(\mathcal{P}_n,r,\rho)$. 
\begin{theorem}\label{Candela12}
Let $A\subset \mathbb{R}^d$ be open and bounded. 
Let $\Sigma$ be one of the soft random simplicial complexes $\Sigma\in\{R(\mathcal{P}_n,r,\rho), C(\mathcal{P}_n,r,\rho)\}$ with multiparameter $\rho=(p_1,p_2,\dots)$. Suppose that $nr^d\to 0$ and $\prod\limits_{i=1}^k p_i^{\binom{k+1}{i+1}}n^{k+1}r^{dk}\to \infty$, then
\begin{align*}
\frac{f_{k,A}^\mathcal{P}(\rho)-E(f_{k,A}^\mathcal{P}(\rho))}{\sqrt{var(f_{k,A}^\mathcal{P}(\rho))}}\xRightarrow{D} N(0,1).   
\end{align*}
\end{theorem}

\begin{proof}[Proof of Theorem~\ref{Candela12}]
We will prove the theorem for $R(\mathcal{P}_n,r,\rho)$. The proof for $C(\mathcal{P}_n,r,\rho)$ is analogous, replacing all instances of $h$ with $\mathfrak{h}$ and  using Theorem~\ref{newcosito} instead of Proposition~\ref{Penrose-1}. We write $f^\mathcal{P}_{i}(\rho,\Sigma)$ simply as $f^\mathcal{P}_{i}(\rho)$ for the remainder of the proof.

Let $\{Q_{m,n}\}_{m\in \mathbb{N}}$ be a partition of $\mathbb{R}^d$ into cubes of side length $r$. Also, if $A\subset \mathbb{R^d}$ is open and bounded, let $I_A$ be the set of indices such that $A\cap Q_{m,n}\neq \emptyset$. Considering the convention from Definition~\ref{h's}, write the number of $k$-faces of $R(\mathcal{P}_n,r,\rho)$ with left most point of the $k$-face in $A$ as
$$f_{k,A}^\mathcal{P}(\rho)=\sum_{m\in I_A}\sum_{Y\subset \mathcal{P}_n}h_{\mathcal{P}_n,r,k+1, A\cap Q_{m,n},\rho}(Y),$$
it follows that, if $$\xi_{m,\rho}=\sum_{Y\subset \mathcal{P}_n}h_{\mathcal{P}_n,r,k+1,(A\cap Q_{m,n}),\rho}(Y),$$
then
$$\frac{f_{k,A}^\mathcal{P}(\rho)-E(f_{k,A}^\mathcal{P}(\rho))}{\sqrt{var(f_{k,A}^\mathcal{P}(\rho))}}=\sum_{m\in I_A} \left (\frac{\xi_{m,\rho}-E(\xi_{m,\rho})}{\sqrt{var(f_{k,A}^\mathcal{P}(\rho))}}\right)=\sum_{m\in I_A}\zeta_m,$$
for 
$$\zeta_m=\frac{\xi_{m,\rho}-E(\xi_{m,\rho})}{\sqrt{var(f_{k,A}^\mathcal{P}(\rho))}}.$$
Define the  relation $\sim$ on $I_A$ by: $m_1\sim m_2$ if and only if the Euclidean distance from $Q_{m_1,n}$ to $Q_{m_2,n}$ is less than $8r$. Thus, $(I_A, \sim)$ is a dependency graph for the collection $\{\zeta_{m,\rho}\}$. Furthermore, the degree of the vertices in this dependency graph is bounded since, given $m_1$, the number of boxes $Q_{m_2,n}$ within a distance $8r$ of $Q_{m_1,n}$ is finite. 

We will use Theorem~\ref{Penrose-4} to guarantee a central limit theorem for $f_{k,A}^\mathcal{P}(\rho)$. For that, we need to prove that, for $q=3,4$, the expression
$$\sqrt{D^{q-1}\sum_{m\in I_A}E|\zeta_m|^q}=\sqrt{D^{q-1}\sum_{m\in I_A}E\left|\frac{\xi_{m,\rho}-E(\xi_{m,\rho})}{\sqrt{var(f_{k,A}^\mathcal{P}(\rho))}}\right|^q}\to 0.$$
Notice that since $A\subset \mathbb{R}^d$ is bounded, the number of $m\in I_A$ (i.e. the indices $m$ such that $A \cap Q_{m,n}\neq \varnothing$) is bounded by a constant times $\frac{1}{r^d}$. Then, for a constant $c$,
$$\sqrt{D^{q-1}\sum_{m\in I_A}E|\zeta_m|^q}=\sqrt{D^{q-1}\sum_{m\in I_A}E\left|\frac{\xi_{m,\rho}-E(\xi_{m,\rho})}{\sqrt{var(f_{k,A}^\mathcal{P}(\rho))}}\right|^q}\leq \sqrt{D^{q-1}\frac{c}{r^d}E\left|\frac{\xi_{m,\rho}-E(\xi_{m,\rho})}{\sqrt{var(f_{k,A}^\mathcal{P}(\rho))}}\right|^q}.$$ 
Thus, in order to use Theorem~\ref{Penrose-4},  it is enough to prove  that, for $q=3,4$,
\begin{equation}\label{b_qn}
b_{q,n}=\frac{1}{r^d}\max_{m}\frac{E|\xi_{m,\rho}- E(\xi_{m,\rho})|^q}{ var(f_{k,A}^\mathcal{P}(\rho))^\frac{q}{2}}\to 0.   
\end{equation}
Notice  that each $\xi_{m,\rho}$ counts the number of $k$-faces of $R(\mathcal{P}_n,r,\rho)$ with left most point of the $k$-face in $A\cap Q_{m,n}$. To show Equation~(\ref{b_qn}),  we need to bound the expressions $|\xi_{m,\rho}- E(\xi_{m,\rho})|^q$. However, by Lemma~\ref{Hyj}, it is enough to bound $E|\xi_{m,\rho}|^q$.
For that, let  $\xi_{m}$ be the number of $k$-faces of $R(\mathcal{P}_n,r)$ with left most point in $A\cap Q_{m,n}$, then by Lemma~\ref{pq}, we have
\begin{equation}\label{99x}
E|\xi_{m,\rho}|^q\leq \prod\limits_{i=1}^k p_i^{\binom{k+1}{i+1}} E|\xi_{m}|^q.    
\end{equation}
Furthermore, if  $z_m$ is the number of points of $\mathcal{P}_n$ within a  distance $2r$ of the cube $Q_{m,n}$, by Lemma~\ref{z}, we have
$$|\xi_{m}|^q\leq (z_m)_{k}$$
and also that, for a constant $c_k$
$$ E|\xi_{m}|^q\leq E((z_m)_k^q)\leq c_k(nr^d)^{k+1}.$$
We now see that
\begin{equation}\label{99x2}
E|\xi_{m,\rho}|^q \leq \prod\limits_{i=1}^k p_i^{\binom{k+1}{i+1}} E|\xi_{m}|^q\leq \prod\limits_{i=1}^k p_i^{\binom{k+1}{i+1}} E((z_m)_k^q)\leq c_k \prod\limits_{i=1}^k p_i^{\binom{k+1}{i+1}}(nr^d)^{k+1}.
\end{equation}
It now follows from Equations (\ref{b_qn}), (\ref{99x}) and (\ref{99x2}) that
$$b_{q,n}\leq \frac{1}{r^d}\max_{m}\frac{\prod\limits_{i=1}^k p_i^{\binom{k+1}{i+1}}(nr^d)^{k+1}}{var(f_{k,A}^\mathcal{P}(\rho))^\frac{q}{2}},$$
and, since by Corollary~\ref{Candela14} we have $var(f_{k,A}^\mathcal{P}(\rho))\sim \mu_{k+1, A}\prod\limits_{i=1}^k p_i^{\binom{k+1}{i+1}}n^{k+1}r^{dk}$, then for some constant $c$ we have
\begin{align}\label{error}
\notag b_{q,n}&\leq \max_{m}\frac{c\prod\limits_{i=1}^k p_i^{\binom{k+1}{i+1}}n^{k+1}r^{dk}}{\left(\prod\limits_{i=1}^k p_i^{\binom{k+1}{i+1}}n^{k+1}r^{dk}\right)^\frac{q}{2}}\\
b_{q,n}&\leq \frac{c}{ \left(\prod\limits_{i=1}^k p_i^{\binom{k+1}{i+1}}n^{k+1}r^{dk}\right)^{\frac{q}{2}-1}}\to 0,
\end{align}
since, by hypothesis, $\prod\limits_{i=1}^k p_i^{\binom{k+1}{i+1}}n^{k+1}r^{dk}\to \infty$. Thus, by Theorem~\ref{Penrose-4}, for $A\subset \mathbb{R^d}$  open and bounded:
\begin{align}\label{clt1}
\frac{f_{k,A}^\mathcal{P}(\rho)-E(f_{k,A}^\mathcal{P}(\rho))}{\sqrt{var(f_{k,A}^\mathcal{P}(\rho))}}\xRightarrow{D} N(0,1).   
\end{align}
\end{proof}
We now move to the case where $A=\mathbb{R}^d$.
\begin{theorem}\label{Candela12xx2}
Let $\Sigma$ be one of the soft random simplicial complexes $\Sigma\in\{R(\mathcal{P}_n,r,\rho), C(\mathcal{P}_n,r,\rho)\}$ with multiparameter $\rho=(p_1,p_2,\dots)$. Suppose that $nr^d\to 0$ and $\prod\limits_{i=1}^k p_i^{\binom{k+1}{i+1}}n^{k+1}r^{dk}\to \infty$, then
$$ \frac{f_k^\mathcal{P}(\rho, \Sigma)-E(f_k^\mathcal{P}(\rho, \Sigma))}{\sqrt{var(f_k^\mathcal{P}(\rho, \Sigma))}}\xRightarrow{D} N(0,1).$$
\end{theorem}
\begin{proof}
We will prove the theorem for $R(\mathcal{P}_n,r,\rho)$. The proof for $C(\mathcal{P}_n,r,\rho)$ is analogous, replacing all instances of $h$ with $\mathfrak{h}$ and  using Theorem~\ref{newcosito} instead of Proposition~\ref{Penrose-1}. We write $f^\mathcal{P}_{i}(\rho,\Sigma)$ simply as $f^\mathcal{P}_{i}(\rho)$ for the remainder of the proof.

Consider the open sets
$$A_K=(-K,K)^d,~~~~A^K=\mathbb{R}^d\backslash [-K,K]^d$$
and the operator $\zeta^{\mathcal{P}}_{n,k}$ which assigns to each subset $A\subset \mathbb{R}^d$ the random variable $\zeta_{n,k}^\mathcal{P}(A)$ defined by
$$\zeta_{n,k}^{\mathcal{P}}(A)=\frac{f_{k,A}^\mathcal{P}(\rho)-E(f_{k,A}^\mathcal{P}(\rho))}{\sqrt{\prod\limits_{i=1}^k p_i^{\binom{k+1}{i+1}}n^{k+1}r^{dk}}}$$
which centralizes a normalization of the random variable $f^\mathcal{P}_{k,A}(\rho)$ around its expected value.
Given $t\in \mathbb{R}$ and $\epsilon>0$, notice that since
$$\{\zeta_{n,k}^{\mathcal{P}}(\mathbb{R}^d)\leq t\}=H_1 \cup H_2\cup H_3$$
for
$$H_1=\{\zeta_{n,k}^{\mathcal{P}}(\mathbb{R}^d)\leq t, |\zeta_{n,k}^{\mathcal{P}}(A_K)-t|<\epsilon\}$$
$$H_2=\{\zeta_{n,k}^{\mathcal{P}}(\mathbb{R}^d)\leq t, \zeta_{n,k}^{\mathcal{P}}(A_K)\geq t+\epsilon\}$$
$$H_3=\{\zeta_{n,k}^{\mathcal{P}}(\mathbb{R}^d)\leq t, \zeta_{n,k}^{\mathcal{P}}(A_K)\leq t-\epsilon\}$$
and 
$$H_3=\{\zeta_{n,k}^{\mathcal{P}}(A_K)\leq t-\epsilon\}\backslash \{\zeta_{n,k}^{\mathcal{P}}(A_K)\leq t-\epsilon, \zeta_{n,k}^{\mathcal{P}}(\mathbb{R}^d)>t\}$$
therefore
\begin{align}\label{000}
P(\zeta_{n,k}^{\mathcal{P}}(\mathbb{R}^d)\leq t)&=P(\zeta_{n,k}^{\mathcal{P}}(A_K)\leq t-\epsilon)- P(\zeta_{n,k}^{\mathcal{P}}(A_K)\leq t-\epsilon, \zeta_{n,k}^{\mathcal{P}}(\mathbb{R}^d)>t)\\
\notag&+P(\zeta_{n,k}^{\mathcal{P}}(\mathbb{R}^d)\leq t, \zeta_{n,k}^{\mathcal{P}}(A_K)\geq t+\epsilon)+ P(\zeta_{n,k}^{\mathcal{P}}(\mathbb{R}^d)\leq t,  |\zeta_{n,k}^{\mathcal{P}}(A_K)-t|<\epsilon).   
\end{align}
On the other hand, the random variable $X$ that counts the number of of points of the Poisson process $\mathcal{P}_n$ in the set $\mathbb{R}^d\backslash (A_K\cup A^K))$  is a Poisson variable of mean $\lambda=\int_{\mathbb{R}^d\backslash (A_K\cup A^K))}nf(x)dx$. Thus, since $Leb(\mathbb{R}^d\backslash (A_K\cup A^K))=0$, by Markov's inequality, Lemma~\ref{Frieze2}, we have
$$P(X>0)\leq E(X)=\int_{\mathbb{R}^d\backslash (A_K\cup A^K))}nf(x)dx=0$$
and therefore, $X=0$ a.a.s. It follows that
$$f_{k,\mathbb{R}^d}^{\mathcal{P}}(\rho)=f_{k,A_K}^{\mathcal{P}}(\rho)+f_{k,A^K}^{\mathcal{P}}(\rho),~~~\text{a.a.s}$$
and
$$\zeta_{n,k}^{\mathcal{P}}(\mathbb{R}^d)= \zeta_{n,k}^{\mathcal{P}}(A_K)+ \zeta_{n,k}^{\mathcal{P}}(A^K),~~~a.a.s.$$
Therefore, a.a.s., we have
\begin{align}\label{001}
P\left(\zeta_{n,k}^{\mathcal{P}}(A_K)\leq t-\epsilon, \zeta_{n,k}^{\mathcal{P}}(\mathbb{R}^d)>t\right)&= P\left(\zeta_{n,k}^{\mathcal{P}}(A_K)\leq t-\epsilon, \zeta_{n,k}^{\mathcal{P}}(A_K)+ \zeta_{n,k}^{\mathcal{P}}(A^K)>t\right)\\
\notag&\leq P\left(\zeta_{n,k}(A^K)>\epsilon\right)  
\end{align}
also,
\begin{align}\label{002}
P\left(\zeta_{n,k}^{\mathcal{P}}(\mathbb{R}^d)\leq t, \zeta_{n,k}^{\mathcal{P}}(A_K)\geq t+\epsilon\right)&=P\left((\zeta_{n,k}^{\mathcal{P}}(A_K)+ \zeta_{n,k}^{\mathcal{P}}(A^K)\leq t, \zeta_{n,k}^{\mathcal{P}}(A_K)\geq t+\epsilon\right)\\
\notag&\leq P\left(\zeta_{n,k}^{\mathcal{P}}(A^K)<-\epsilon\right)   
\end{align}
and
\begin{align}\label{003}
P\left(\zeta_{n,k}^{\mathcal{P}}(\mathbb{R}^d)\leq t, |\zeta_{n,k}^{\mathcal{P}}(A_K)-t|<\epsilon\right)&=P\left(\zeta_{n,k}^{\mathcal{P}}(A_K)+ \zeta_{n,k}^{\mathcal{P}}(A^K)\leq t, |\zeta_{n,k}^{\mathcal{P}}(A_K)-t|<\epsilon\right))\\
\notag&\leq P\left(|\zeta_{n,k}^{\mathcal{P}}(A_K)-t|<\epsilon\right).
\end{align}
Therefore, from Equation (\ref{000}), we have
\begin{align}\label{der}
|P\left(\zeta_{n,k}^{\mathcal{P}}(\mathbb{R}^d)\leq t\right)- P\left(\zeta_{n,k}^{\mathcal{P}}(A_K)\leq t-\epsilon\right)|\leq P\left(|\zeta_{n,k}^{\mathcal{P}}(A_K)-t|<\epsilon\right) + P\left(|\zeta_{n,k}^{\mathcal{P}}(A^K)|>\epsilon\right),    
\end{align}
where the last term combines Equations (\ref{002}) and (\ref{003}).
By Chebychev's inequality, Lemma~\ref{Frieze3}, and the central limit theorem established for bounded sets in Equation (\ref{clt1}), the right hand of Equation (\ref{der}) is bounded above by 
\begin{align}\label{der2}
\frac{1}{\epsilon^2}var(\zeta_{n,k}^{\mathcal{P}}(A^K))+ P\left(\left|\sqrt{\frac{var(f_{k,A_K}^\mathcal{P}(\rho))}{\prod\limits_{i=1}^k p_i^{\binom{k+1}{i+1}}n^{k+1}r^{dk}}}Z-t\right|<\epsilon\right)+\frac{c}{ \left(\prod\limits_{i=1}^k p_i^{\binom{k+1}{i+1}}n^{k+1}r^{dk}\right)^{\frac{q}{2}-1}}   
\end{align}
where $Z$ is a normal standard variable. The second and third expressions in Equation~(\ref{der2}) come from the normal approximation theorem stated in (\ref{clt1}) and the error for using it, Equation~(\ref{error}). Notice now , since $Z$ is a standard normal variable, the second term in Equation~(\ref{der2}) becomes
\begin{align}\label{der3}
P\left(\left|\sqrt{\frac{var(f_{k,A_K}^\mathcal{P}(\rho))}{\prod\limits_{i=1}^k p_i^{\binom{k+1}{i+1}}n^{k+1}r^{dk}}}Z-t\right|<\epsilon\right)&=P\left((t- \epsilon)\sqrt{\frac{ \prod\limits_{i=1}^{k}p_i^{\binom{k+1}{i+1}}n^{k+1}r^{dk}}{var(f_{k,A_K}^\mathcal{P}(\rho))}} \leq Z\leq (t+ \epsilon)\sqrt{\frac{\prod\limits_{i=1}^{k}p_i^{\binom{k+1}{i+1}}n^{k+1}r^{dk}}{var(f_{k,A_K}^\mathcal{P}(\rho))}} \right)\\
\notag&\leq 2\epsilon\sqrt{\frac{n^{k+1}r^{dk}}{2\pi \cdot var(f_{k,A_K}^\mathcal{P}(\rho))}}.
\end{align}
Furthermore, by Corollary~\ref{Candela14}, $$var(f_{k,A_K}^\mathcal{P}(\rho))=\mu_{k+1, A_K}\prod\limits_{i=1}^k p_i^{\binom{k+1}{i+1}}n^{k+1}r^{dk}$$
and 
$$var(\zeta_{n,k}(A_K)^{\mathcal{P}})=\mu_{k+1, A^K}.$$ 
Thus, from Equations (\ref{der}),(\ref{der2}) and (\ref{der3}), for $K$ and $\epsilon$ fixed we have 
\begin{align}\label{yeah}
\limsup_{n\to\infty}|P(\zeta_{n,k}^{\mathcal{P}}(\mathbb{R}^d)\leq t)- P(\zeta_{n,k}^{\mathcal{P}}(A_K)\leq t-\epsilon)|\leq \frac{1}{\epsilon^2}\mu_{k+1, A^K}+ 2\epsilon\sqrt{\frac{1}{2\pi \mu_{k+1, A_K}}}.   
\end{align}
On the other hand, using the central limit theorem established for bounded sets in Equation (\ref{clt1}), Equation (\ref{yeah}) is equivalent to
$$\limsup_{n\to\infty}\left|P(\zeta_{n,k}^{\mathcal{P}}(\mathbb{R}^d)\leq t)- P\left(\sqrt{\frac{var(f_{k,A_K}^\mathcal{P}(\rho))}{\prod\limits_{i=1}^k p_i^{\binom{k+1}{i+1}}n^{k+1}r^{dk}}}Z\leq t-\epsilon\right)\right|\leq \frac{1}{\epsilon^2}\mu_{k+1, A^K}+ 2\epsilon\sqrt{\frac{1}{2\pi \mu_{k+1, A_K}}}.$$
Now, if $\Phi$ is the cumulative distribution function of the normal standard variable $Z$ then
$$P\left(\sqrt{\frac{var(f_{k,A_K}^\mathcal{P}(\rho)))}{\prod\limits_{i=1}^k p_i^{\binom{k+1}{i+1}}n^{k+1}r^{dk}}}Z\leq t-\epsilon\right)= \Phi\left((t-\epsilon)\sqrt{\frac{\prod\limits_{i=1}^k p_i^{\binom{k+1}{i+1}}n^{k+1}r^{dk}}{var(f_{k,A_K}^\mathcal{P}(\rho))}}\right)\to \Phi\left((t-\epsilon)\sqrt{\frac{1}{\mu_{k+1, A_K}}}\right).$$
Therefore,
$$\limsup_{n\to\infty}\left|P(\zeta_{n,k}^{\mathcal{P}}(\mathbb{R}^d)\leq t)- \Phi\left((t-\epsilon)\sqrt{\frac{1}{\mu_{k+1, A_K}}}\right)\right|\leq  \frac{1}{\epsilon^2}\mu_{k+1,A^K}+ 2\epsilon\sqrt{\frac{1}{2\pi \mu_{k+1,A_K}}}.$$
Recall that, by Equation (\ref{newmu}),
$$\mu_{k+1, A}=\frac{1}{(k+1)!}\int_A f(x)^kdx \int_{(\mathbb{R}^d)^{k}}h_{n,1,\Gamma}(\{0,x_1,...,x_{k}\})~dx_1dx_2...dx_{k},$$
and therefore, by definition of $A_K$ and $A^K$,  $\lim\limits_{K\to \infty}\mu_{k+1,A_K}=0$ and $\lim\limits_{K\to \infty}\mu_{k+1, A^K}=\mu_{k+1}$, where $\mu_{k+1}$ stands for $\mu_{k+1,\mathbb{R}^d}$ . Thus, for $K$ large enough and $\epsilon'>0$, we have 
\begin{align}\label{huhu2}
\limsup_{n\to\infty}\left|P(\zeta_{n,k}^{\mathcal{P}}(\mathbb{R}^d)\leq t)- \Phi\left(\sqrt{(t-\epsilon)\frac{1}{\mu_{k+1, A_K}}}\right)\right|\leq \epsilon+ 2\epsilon\sqrt{\frac{1}{2\pi (\mu_{k+1}-\epsilon')}}.
\end{align}
Finally, consider the expression 
\begin{align*}
&\quad\limsup_{n\to\infty}\left|P(\zeta_{n,k}^{\mathcal{P}}(\mathbb{R}^d)\leq t)- \Phi\left(\sqrt{(t-\epsilon)\frac{1}{\mu_{k+1, A_K}}}\right)\right|\\
&=\limsup_{n\to\infty}\left|P(\zeta_{n,k}^{\mathcal{P}}(\mathbb{R}^d)\leq t)- \Phi\left(t\sqrt{\frac{1}{\mu_{k+1}}}\right)+\Phi\left(t\sqrt{\frac{1}{\mu_{k+1}}}\right)- \Phi\left(\sqrt{(t-\epsilon)\frac{1}{\mu_{k+1, A_K}}}\right)\right|\\
&\geq\limsup_{n\to\infty}\left|\left|P(\zeta_{n,k}^{\mathcal{P}}(\mathbb{R}^d)\leq t)- \Phi\left(t\sqrt{\frac{1}{\mu_{k+1}}}\right)\right|-\left| \Phi\left(\sqrt{(t-\epsilon)\frac{1}{\mu_{k+1, A_K}}}\right)-\Phi\left(t\sqrt{\frac{1}{\mu_{k+1}}}\right)\right|\right|.
\end{align*}
Since as $\epsilon\to 0$ and $K\to \infty$, $\Phi\left((t-\epsilon)\sqrt{\frac{1}{\mu_{k+1,A_K}}}\right)\to \Phi\left(t\sqrt{\frac{1}{\mu_{k+1}}}\right)$ then, for $\epsilon_1>0$, we have that 
\begin{align}\label{huhu}
\limsup_{n\to\infty}\left|P(\zeta_{n,k}^{\mathcal{P}}(\mathbb{R}^d)\leq t)- \Phi\left(\sqrt{(t-\epsilon)\frac{1}{\mu_{k+1, A_K}}}\right)\right|\\
\notag\geq\limsup_{n\to\infty}\left|\left|P(\zeta_{n,k}^{\mathcal{P}}(\mathbb{R}^d)\leq t)- \Phi\left(t\sqrt{\frac{1}{\mu_{k+1}}}\right)\right|-\epsilon_1\right|.    
\end{align}
Thus, by equations (\ref{huhu2}) and (\ref{huhu}) combined we have  
$$\limsup_{n\to\infty}\left|\left|P(\zeta_{n,k}^{\mathcal{P}}(\mathbb{R}^d)\leq t)- \Phi\left(t\sqrt{\frac{1}{\mu_{k+1}}}\right)\right|-\epsilon_1\right|\leq \epsilon+ 2\epsilon\sqrt{\frac{1}{2\pi (\mu_{k+1}-\epsilon')}},$$
and since $\epsilon$, $\epsilon'$ and $\epsilon_1$ were arbitrary, it follows that
$$\limsup_{n\to\infty}\left|P(\zeta_{n,k}^{\mathcal{P}}(\mathbb{R}^d)\leq t)- \Phi\left(t\sqrt{\frac{1}{\mu_{k+1}}}\right)\right|=0,$$
and then
$$\lim_{n\to\infty}\left|P(\zeta_{n,k}^{\mathcal{P}}(\mathbb{R}^d)\leq t)- \Phi\left(t\sqrt{\frac{1}{\mu_{k+1}}}\right)\right|=0.$$
Therefore, 
$$\zeta_{n,k}^{\mathcal{P}}(\mathbb{R}^d)\xRightarrow{D} \sqrt{\mu_{k+1}}Z,$$ 
and this is
$$\frac{f_{k}^\mathcal{P}(\rho)-E(f_{k}^\mathcal{P}(\rho))}{\sqrt{var(f_{k}^\mathcal{P}(\rho))}}\xRightarrow{D} N(0,1).$$
\end{proof}

\subsection{De-Poissonization}
Finally, the remaining work is to use a de-Poissonization argument to guarantee that the central limit theorem proved above for $f^\mathcal{P}_k(\rho)$ is also true for $f_k(\rho)$ itself. For that, we use the following result:
\begin{theorem}[Theorem 2.12 \cite{MR1986198}]\label{Penrose-6}
Suppose that for each $n\in \mathbb{N}$, $H_n(\mathcal{X})$ is a real-valued functional on finite sets $\mathcal{X}\subset \mathbb{R}^d$. Suppose that for some $\sigma^2>0$,
\begin{enumerate}
\item $\frac{1}{n} var(H_n(\mathcal{P}_n))\to \sigma^2, and$
\item $\frac{1}{\sqrt{n}} (H_n(\mathcal{P}_n- E(H_n(\mathcal{P}_n))\to \sigma^2Z$, for $Z$ a standard normal random variable. 
\end{enumerate}
Suppose that there exist constants $\alpha\in \mathbb{R}$ and $\gamma>\frac{1}{2}$ such that the increments $R_{m,n}=H_n(\mathcal{X}_{m+1})-H_n(\mathcal{X}_m)$ satisfy
$$\lim_{n\to\infty}\left( \sup_{n-n^\gamma\leq m\leq n+n^\gamma}|E(R_{m,n})-\alpha|)\right) =0$$
$$\lim_{n\to\infty}\left( \sup_{n-n^\gamma\leq m<m'\leq n+n^\gamma}|E(R_{m,n}R_{m,n})-\alpha^2|)\right) =0$$
and
$$\lim_{n\to\infty}\left( \frac{1}{\sqrt{n}}\sup_{n-n^\gamma\leq m\leq n+n^\gamma}E(R_{m,n}^2)\right)=0.$$
Finally, assume that there is a  constant $\beta>0$ such that, with probability one, for $\mathcal{X}_n=\{X_1,\dots,X_n\}$, we have
$$|H_n(\mathcal{X}_{m})|\leq \beta(n+m)^\beta.$$
Then $\alpha^2\leq \sigma^2$, and as $n\to \infty$, $\frac{1}{n}var(H_n(\mathcal{X}_{n}))\to \sigma^2-\alpha^2$ and
$$\frac{1}{\sqrt{n}}(H_n(\mathcal{X}_{n}- E(H_n(\mathcal{X}_{n}))\to Z\sqrt{\sigma^2-\alpha^2}.$$
\end{theorem}

\begin{proof}[\textbf{Proof of Theorem~\ref{Candela15}} ]
We will prove the theorem for $R(\mathcal{P}_n,r,\rho)$. The proof for $C(\mathcal{P}_n,r,\rho)$ is analogous, replacing all instances of $h$ with $\mathfrak{h}$ and  using Theorem~\ref{newcosito} instead of Proposition~\ref{Penrose-1}. We write $f_{i}(\rho,\Sigma)$ simply as $f_{i}(\rho)$ for the remainder of the proof.

Using the conventions from Definition \ref{h's}, for $Y\subset \mathcal{X}_n$ with $|Y|=k+1$, define  the following functional
\begin{align}\label{aguanta}
H_n( \mathcal{X}_n)=\sqrt{n}\frac{\sum\limits_{Y\subset \mathcal{X}_n}h_{n,r,k+1,\rho}(Y)}{\sqrt{\prod\limits_{i=1}^kp_i^{\binom{k+1}{i+1}}n^{k+1}r^{dk}}}.    
\end{align}
Consider the increments $R_{m,n}$ of the functional $H$ given by 
\begin{align*}
R_{m,n}&=H_n(\mathcal{X}_{m+1})- H_n(\mathcal{X}_{m})\\
&=\frac{\sqrt{n}}{\sqrt{\prod\limits_{i=1}^kp_i^{\binom{k+1}{i+1}}n^{k+1}r^{dk}}}\left(\sum\limits_{Y\subset \mathcal{X}_{m+1}}h_{n,r,k+1,\rho}(Y)- \sum\limits_{Y\subset \mathcal{X}_m}h_{n,r,k+1,\rho}(Y)\right)\\
&=\frac{1}{\sqrt{\prod\limits_{i=1}^kp_i^{\binom{k+1}{i+1}}n^{k}r^{dk}}}D_{m,n}^k
\end{align*}
where
$$D_{m,n}^k=\sum\limits_{Y\subset \mathcal{X}_{m+1}}h_{n,r,k+1,\rho}(Y)- \sum\limits_{Y\subset \mathcal{X}_m}h_{n,r,k+1,\rho}(Y).$$
Note that $D_{m,n}^k$ counts the number of $k$-faces in $R(\mathcal{X}_{m+1},r,\rho)$  that contain the vertex $X_{m+1}$. Thus, by Theorem~  \ref{hvsh}:
\begin{align}\label{sofia2}
E(D_{m,n}^k)=\binom{m}{k}E(h_{n,r,k+1,\rho}(\mathcal{X}_{k+1}))= \binom{m}{k}E(h_{n,r,k+1}(\mathcal{X}_{k+1}))\prod\limits_{i=1}^kp_i^{\binom{k+1}{i+1}}    
\end{align}
thus, for $n-n^{\frac{2}{3}}\leq  m \leq n+n^{\frac{2}{3}}$,
\begin{align*}
1=\lim_{n\to \infty} \frac{n-n^{\frac{2}{3}}}{n}\leq \lim_{n\to \infty} \frac{m}{n} \leq \lim_{n\to \infty} \frac{n+n^{\frac{2}{3}}}{n}=1
\end{align*}
and therefore, $m\sim n$. Thus
\begin{align}\label{sofia}
\binom{m}{k}=\frac{m!}{(m-k)!k!}\sim \frac{n^k}{k!}.   
\end{align}
Furthermore, by Proposition~\ref{Penrose-1} we have
\begin{align}\label{bo}
E(h_{n,r,k+1}(\mathcal{X}_{k+1}))\sim (k+1)! \mu_{k+1}r^{dk}.    
\end{align}
Therefore, by Equations~(\ref{sofia2}), (\ref{sofia}) and (\ref{bo}) combined we have
$$\binom{m}{k}E(h_{n,r,k+1}(\mathcal{X}_{k+1}))\prod\limits_{i=1}^kp_i^{\binom{k+1}{i+1}}\sim
\frac{n^k}{k!}(k+1)! \prod\limits_{i=1}^kp_i^{\binom{k+1}{i+1}}\mu_{k+1}r^{dk},$$
It follows that, uniformly over $n-n^{\frac{2}{3}}\leq  m \leq n+n^{\frac{2}{3}}$, we have
\begin{equation}\label{E1}
E(D_{m,n}^k)=\binom{m}{k}E(h_{n,r,k+1,\rho}(\mathcal{X}_{k+1}))\sim (k+1)\prod\limits_{i=1}^kp_i^{\binom{k+1}{i+1}}(nr^d)^{k}\mu_{k+1}.
\end{equation}
Therefore,by definition of $R_{m,n}$,  we have
\begin{align*}
E(R_{m,n})&= \frac{1}{\sqrt{\prod\limits_{i=1}^kp_i^{\binom{k+1}{i+1}}n^{k}r^{dk}}} E(D_{m,n}^k)\\
&\sim \frac{ (k+1)\prod\limits_{i=1}^kp_i^{\binom{k+1}{i+1}}(nr^d)^{k}\mu_{k+1}}{\sqrt{\prod\limits_{i=1}^kp_i^{\binom{k+1}{i+1}}n^{k}r^{dk}}}\\
&\sim  (k+1)\mu_{k+1} \sqrt{\prod\limits_{i=1}^kp_i^{\binom{k+1}{i+1}}(nr^d)^{k}}\to 0,
\end{align*}
since by hypothesis $nr^d\to 0$. Therefore, for $\alpha=0$ and $\gamma=\frac{2}{3}$, we have
\begin{equation}\label{lim1}
\lim_{n\to\infty}\left( \sup_{n-n^\gamma\leq m\leq n+n^\gamma}|E(R_{m,n})-\alpha|)\right)\to 0.
\end{equation}
We now examine $E(D_{m,n}^kD_{l,n}^k)$. Recall from the above that  that $D_{m,n}^k$ counts the number of $k$-faces in $R(\mathcal{X}_{m+1},r,\rho)$  that contain the vertex $X_{m+1}$.
Thus, if $m<l$, then
\begin{equation}\label{dobled's}
E(D_{m,n}^kD_{l,n}^k)= E\left(\sum_{Y\subset \mathcal{X}_{m} }\sum_{Y'\subset \mathcal{X}_{l} }h_{n,r,k+1,\rho}(Y\cup\{X_{m+1}\})h_{n,r,k+1,\rho}(Y'\cup\{X_{l+1}\})\right)
\end{equation}
$$= E\left(\sum_{j=0}^{k}\sum_{Y\subset \mathcal{X}_{m} }\sum_{Y'\subset \mathcal{X}_{l} }h_{n,r,k+1,\rho}(Y\cup\{X_{m+1}\})h_{n,r,k+1,\rho}(Y'\cup\{X_{l+1}\})1_{\{|(Y\cup\{X_{m+1})\cap Y'\}|=j\}}\right).$$
$$=E\left(\sum_{j=0}^{k}\binom{m}{k}\binom{k+1}{j}\binom{l-k-1}{k+1-j-1}h_{n,r,k+1,\rho}(X_{k+1})h_{n,r,k+1,\rho}(X'_{k+1})1_{\{|(X'_{k+1}\cap X_{k+1}\}|=j\}}\right)$$
Let us consider first the term $j=0$ in the sum above. This is given by the expression
$$\alpha_0=\binom{m}{k}\binom{l-k-1}{k}E\left(h_{n,r,k+1,\rho}(Y\cup\{X_{m+1}\})h_{n,r,k+1,\rho}(Y'\cup\{X_{l+1}\})1_{\{|(Y\cup\{X_{m+1})\cap Y'\}|=0\}}\right).$$
Notice now that since $j=0$, $(Y\cup\{X_{m+1})\cap Y'=\emptyset$, and therefore the functions $h_{n,r,k+1,\rho}(Y\cup\{X_{m+1}\})$ and $h_{n,r,k+1,\rho}(Y'\cup\{X_{l+1}\})$ are independent. it follows that 
$$\alpha_0=\binom{m}{k}\binom{l-k-1}{k}E(\left(h_{n,r,k+1,\rho}(Y\cup\{X_{m+1}\}))E(h_{n,r,k+1,\rho}(Y'\cup\{X_{l+1}\})\right),$$
therefore, for $m,l\in [n-n^{\frac{2}{3}}, n+n^{\frac{2}{3}}]$, by Equation (\ref{E1})
\begin{align}\label{alpha0}
\alpha_0\sim(k+1)^2\prod\limits_{i=1}^kp_i^{2\binom{k+1}{i+1}}(nr^d)^{2k}\mu^2_{k+1}.
\end{align}
For $j\geq 1$, notice that the following sum $S$ satisfies
\begin{align*}
S&\coloneqq E\left(\sum_{j=1}^{k}\sum_{Y\subset \mathcal{X}_{m} }\sum_{Y'\subset \mathcal{X}_{l} }h_{n,r,k+1,\rho}(Y\cup\{X_{m+1}\})h_{n,r,k+1,\rho}(Y'\cup\{X_{l+1}\})1_{\{|(Y\cup\{X_{m+1})\cap Y'\}|=j\}}\right)\\
&=(\sum_{j=1}^{k}\binom{m}{k}\binom{k+1}{j}\binom{l-k-1}{k+1-j-1}h_{n,r,k+1,\rho}(\mathcal{X}_{k}\cup \{X_m+1\})\\
&\times h_{n,r,k+1,\rho}(\mathcal{X}'_{k}\cup \{X_{l+1}\})1_{\{|(\mathcal{X}_k\cup \{X_{m+1}\})\cap \mathcal{X}'_{k}|=j\}}).   
\end{align*}
Therefore, for $m,l\in [n-n^{\frac{2}{3}}, n+n^{\frac{2}{3}}]$ there exists a constant $c$ such that
$$S\sim \sum_{j=1}^{k} c n^{2k-j}E\left(h_{n,r,k+1,\rho}(\mathcal{X}_k\cup \{X_{m+1}\})h_{n,r,k+1,\rho}(\mathcal{X}'_k\cup \{X_{l+1}\})1_{\{|(\mathcal{X}_k\cup \{X_{m+1}\})\cap \mathcal{X}'_{k}|=j\}}\right).$$
It follows from Theorem~\ref{prodh's} that
\begin{align}\label{ssssss}
S\sim \sum_{j=1}^{k} c n^{2k-j}\prod\limits_{i=1}^kp_i^{2\binom{k+1}{i+1}-\binom{j}{i+1}}E\left(h_{n,r,k+1}(\mathcal{X}_k\cup \{X_{m+1}\})h_{n,r,k+1}(\mathcal{X}'_k\cup \{X_{l+1})1_{\{|(\mathcal{X}_k\cup \{X_{m+1}\})\cap \mathcal{X}'_{k}|=j\}}\right).  
\end{align}
Note that the configurations $\mathcal{X}_k\cup \{X_{m+1}\}$ and $\mathcal{X}'_k\cup \{X_{l+1}\}$ which contribute to the expected value in Equation (\ref{ssssss}) above must form a connected graph of $2k+2-j$ vertices. Combining this with Proposition~\ref{Penrose-1}, we have 
$$S\sim \sum_{j=1}^{k} c n^{2k-j}\prod\limits_{i=1}^kp_i^{2\binom{k+1}{i+1}-\binom{j}{i+1}} r^{d(2k+2-j-1)}.$$
Notice that in the above expression the the final term dominates the rest, then $S$ is asymptotically equal to 
\begin{align}\label{alpha}
S \sim  c n^{k}\prod\limits_{i=1}^kp_i^{2\binom{k+1}{i+1}-\binom{k}{i+1}} r^{d(k-1)}. 
\end{align}
Finally, let us compare the term $\alpha_0$ ($j=0$) with the asymptotic expression (Equation (\ref{alpha})) we found for the sum of the rest of the $j$'s. By equations (\ref{alpha0}) and (\ref{alpha}),
\begin{align*}
\frac{S}{\alpha_0}&\sim \frac{c n^{k}\prod\limits_{i=1}^kp_i^{2\binom{k+1}{i+1}-\binom{k}{i+1}} r^{d(k-1)}}{(k+1)^2\prod\limits_{i=1}^kp_i^{2\binom{k+1}{i+1}}(nr^d)^{2k}\mu^2_{k+1}}\\
&\leq \frac{c n^{k}\prod\limits_{i=1}^kp_i^{\binom{k+1}{i+1}} r^{d(k-1)}}{(k+1)^2\prod\limits_{i=1}^kp_i^{2\binom{k+1}{i+1}}(nr^d)^{2k}\mu^2_{k+1}}\\
&\leq \frac{cnr^d }{(k+1)^2\prod\limits_{i=1}^kp_i^{\binom{k+1}{i+1}}n^{k+1}r^{dk}\mu^2_{k+1}} 
\end{align*}
and since by hypothesis, $nr^d\to 0$ and $\prod\limits_{i=1}^kp_i^{\binom{k+1}{i+1}}n^{k+1} r^{dk}\mu^2_{k+1}\to \infty$, then we have
$$\frac{S}{\alpha_0}\leq\frac{cnr^d }{(k+1)^2\prod\limits_{i=1}^kp_i^{\binom{k+1}{i+1}}n^{k+1} r^{dk}\mu^2_{k+1}}\to 0.$$
Therefore, uniformly over $m,l\in [n-n^{\frac{2}{3}}, n+n^{\frac{2}{3}}]$, as $n\to \infty,$
\begin{equation}\label{proD's}
E(D_{m,n}^kD_{l,n}^k)\sim (k+1)^2\prod\limits_{i=1}^kp_i^{2\binom{k+1}{i+1}}(nr^d)^{2k}\mu_{k+1}^2.
\end{equation}
Thus by hypothesis 
\begin{align*}
 E(R_{m,n}R_{l,n})&= \frac{1}{\prod\limits_{i=1}^kp_i^{\binom{k+1}{i+1}}n^{k}r^{dk}} E(D_{m,n}^kD_{l,n}^k)\\
 &\sim \frac{ (k+1)^2\prod\limits_{i=1}^kp_i^{2\binom{k+1}{i+1}}(nr^d)^{2k}\mu_{k+1}^2}{\prod\limits_{i=1}^kp_i^{\binom{k+1}{i+1}}n^{k}r^{dk}} \\
 &\sim (k+1)^2 \mu_{k+1}^2\prod\limits_{i=1}^kp_i^{\binom{k+1}{i+1}}(nr^d)^{k}\to 0.
\end{align*}
Therefore, for $\alpha=0$ and $\gamma=\frac{2}{3}$,
\begin{equation}\label{lim2}
\lim_{n\to\infty}\left( \sup_{n-n^\gamma\leq m<l\leq n+n^\gamma}|E(R_{m,n}R_{l,n})-\alpha^2|)\right)\to 0.
\end{equation}
Furthermore, uniformly over $m\in [n-n^{\frac{2}{3}},  n+n^{\frac{2}{3}}]$, as $n\to \infty$, 
\begin{align*}
\frac{1}{\sqrt{n}}E(R_{m,n}^2)& \sim (k+1)^2 \mu_{k+1}^2\prod\limits_{i=1}^kp_i^{\binom{k+1}{i+1}}\frac{(nr^d)^{k}}{\sqrt{n}}\\
&\sim (k+1)^2 \mu_{k+1}^2\prod\limits_{i=1}^kp_i^{\binom{k+1}{i+1}}n^{k-\frac{1}{2}}r^{kd}\\
&\leq  (k+1)^2 \mu_{k+1}^2\prod\limits_{i=1}^kp_i^{\binom{k+1}{i+1}}n^{k-\frac{1}{2}}r^{d(k-\frac{1}{2})}\to 0. 
\end{align*}
Therefore, since $r \to 0$ as $n \to \infty$,
\begin{equation}\label{lim3}
\lim_{n\to\infty}\left( \frac{1}{\sqrt{n}}\sup_{n-n^\gamma\leq m\leq n+n^\gamma}|E(R_{m,n}^2|)\right)\to 0.
\end{equation}
By the estimates in Equations (\ref{lim1}), (\ref{lim2}), and (\ref{lim3}), together with Theorem~\ref{Candela12xx2}, the functional $H$ defined in Equation~\ref{aguanta} satisfies all the conditions of  Theorem~\ref{Penrose-6} for $\alpha=0$. Therefore, by Theorem~\ref{Penrose-6},%
for $\sigma$ such that
$$\frac{var(H_n( \mathcal{P}_n))}{n}\to \sigma^2,$$
we have
$$\frac{H_n( \mathcal{X}_n)- E(H_n( \mathcal{X}_n))}{\sqrt{n}}\Rightarrow N(0, \sigma^2).$$
That is,
$$\frac{ \sum\limits_{Y\subset \mathcal{X}_n}\left(h_{n,r,k+1,\rho}(Y)-E(h_{n,r,k+1,\rho}(Y))\right)}{\sqrt{\prod\limits_{i=1}^kp_i^{\binom{k+1}{i+1}}n^{k+1}r^{dk}}}\Rightarrow N(0, \sigma^2).$$
Thus
$$\frac{ f_k(\rho)-E(f_k(\rho))}{\sqrt{\prod\limits_{i=1}^kp_i^{\binom{k+1}{i+1}}n^{k+1}r^{dk}}}\Rightarrow N(0, \sigma^2).$$
Finally, by Theorem~\ref{Candela12xx2}, we have
$$\frac{var(H_n( \mathcal{P}_n))}{n}\sim var(f_k^\mathcal{P}-E(f_k^\mathcal{P}))\sim \mu_{k+1}=\sigma^2.$$
Therefore, 
$$\frac{ f_k(\rho)-E(f_k(\rho))}{\sqrt{\mu_{k+1}\prod\limits_{i=1}^kp_i^{\binom{k+1}{i+1}}n^{k+1}r^{dk}}}\xRightarrow{D} N(0, 1),$$
as desired.
\end{proof}

\section{C.L.T for the Euler characteristic of $R(\mathcal{P}_n,r,\rho)$ and $C(\mathcal{P}_n,r,\rho)$}\label{section5}

In this section, we prove Theorem~\ref{Candela16}, which establishes a central limit theorem  for the Euler characteristic of the soft random simplicial complexes $R(n,r,\rho)$ and $C(n,r,\rho)$. As before, we begin with the Poisson version of the problem for bounded sets.

\begin{theorem}\label{Candela9}
Let $A\subset\mathbb{R}^d$ be  open and bounded. Let $\Sigma$ be one of the soft random simplicial complexes $\Sigma\in\{R(\mathcal{P}_n,r,\rho), C(\mathcal{P}_n,r,\rho)\}$ with multiparameter $\rho=(p_1,p_2,\dots)$, and  let $\chi_A(\rho,\Sigma)$ be defined as
$$\chi_A^\mathcal{P}(\rho, \Sigma)=\sum_{i=0}^{n-1}(-1)^{i+1}f_{i,A}^\mathcal{P}(\rho,\Sigma).$$
Assume that $nr^d\to 0$ and that there exists a non-negative integer $l$ such that
$$\prod\limits_{i=1}^{l+1} p_i^{\binom{l+2}{i+1}}n^{l+2}r^{d(l+1)}\to 0~~  and~~\prod\limits_{i=1}^l p_i^{\binom{l+1}{i+1}}n^{l+1}r^{dl}\to \infty.$$
Then 
$$ \frac{\chi_A^\mathcal{P}(\rho, \Sigma)-E(\chi_A^\mathcal{P}(\rho, \Sigma))}{\sqrt{var(\chi_A^\mathcal{P}(\rho, \Sigma))}}\xRightarrow{D} N(0,1).$$
\end{theorem}

\begin{proof}
We will prove the theorem for $R(\mathcal{P}_n,r,\rho)$ (the proof for $C(\mathcal{P}_n,r,\rho)$ is analogous, replacing all instances of $h$ with $\mathfrak{h}$) and  using Theorem~\ref{newcosito} instead of Proposition~\ref{Penrose-1}. 
We write $\chi^\mathcal{P}_{A}(\rho,\Sigma)$ and $f^\mathcal{P}_{i,A}(\rho,\Sigma)$ simply as $\chi^\mathcal{P}_{A}(\rho)$ and  $f^\mathcal{P}_{i,A}(\rho)$, respectively, for the remainder of the proof.

By Lemma~\ref{Candela11},  under the hypothesis of  Theorem~\ref{Candela9}, a.a.s we have
$$\chi_A^\mathcal{P}(\rho)=f_{0,A}^\mathcal{P}(\rho)-f_{1,A}^\mathcal{P}(\rho)+f_{2,A}^\mathcal{P}(\rho)-f_{3,A}^\mathcal{P}(\rho,)+...+(-1)^{l+1}f_{l,A}^\mathcal{P}(\rho)=\sum_{i=0}^l(-1)^{i+1}f_{i,A}^\mathcal{P}(\rho).$$
Let $\{Q_{m,n}\}_{m\in \mathbb{N}}$ be a partition of $\mathbb{R}^d$ into cubes of side length $r$. Also, if $A\subset \mathbb{R^d}$ is open and bounded,  let $I_A$ be the set of indices such that $A\cap Q_{m,n}\neq \emptyset$. Considering the convention from Definition~\ref{h's}, write the number of $k$-faces of  $R(\mathcal{P}_n,r,\rho)$ with left most point of the $k$-face in $A$ as
$$f_{k,A}^\mathcal{P}(\rho)=\sum_{m\in I_A}\sum_{Y\subset \mathcal{P}_n}h_{\mathcal{P}_n,r,k+1,(A\cap Q_{m,n}),\rho}(Y).$$
It follows that, defining
$$\xi_{m,j,\rho}\coloneqq\sum_{Y\subset \mathcal{P}_n}(h_{\mathcal{P}_n,r,j+1,(A\cap Q_{m,n}),\rho}(Y),$$
we have
$$\frac{\chi_A^\mathcal{P}(\rho)-E(\chi_A^\mathcal{P}(\rho))}{\sqrt{var(\chi_A^\mathcal{P}(\rho))}}=\sum_{m\in I_A} \sum_{j=0}^{l}(-1)^{j+1}\left(\frac{\xi_{m,j,\rho}-E(\xi_{m,j,\rho})}{\sqrt{var(\chi_A^\mathcal{P}(\rho))}}\right)= \sum_{m\in I_A} \frac{\eta_{m,\rho}-E(\eta_{m,\rho})}{\sqrt{var(\chi_A^\mathcal{P}(\rho))}}=\sum_{m\in I_A} \zeta_{m}$$
for 
$$\eta_{m,\rho}\coloneqq \sum_{j=0}^{l}(-1)^{j+1}(\xi_{m,j,\rho}).$$
and 
$$\zeta_{m}\coloneqq \frac{\eta_{m,\rho}-E(\eta_{m,\rho})}{\sqrt{var(\chi_A^\mathcal{P}(\rho))}}.$$
Define the  relation $\sim$ on $I_A$ by: $m_1\sim m_2$ if and only if the Euclidean distance from $Q_{m_1,n}$ to $Q_{m_2,n}$ is less than $8r$. Thus $(I_A, \sim)$ is a dependency graph (Definition~\ref{dependency}) for the collection $\{\zeta_m\}$. Furthermore, the degree of the vertices in this dependency graph is bounded by a constant since, given $m_1$, the number of boxes $Q_{m_2,n}$ within a distance $8r$ of $Q_{m_1,n}$ is finite.

We will use Theorem~\ref{Penrose-4} to derive a central limit theorem for $\chi_A^\mathcal{P}(\rho)$. Notice that since $A\subset \mathbb{R}^d$ is bounded, by Equation~(\ref{b_qn}), it is enough to prove  that, for $q=3,4$,
\begin{align}\label{bq2}
b_{q,n}=\frac{1}{r^d}\max_{m}\frac{E|\eta_{m,\rho}-E(\eta_{m,\rho}))|^q}{ var(\chi_A^\mathcal{P}(\rho))^\frac{q}{2}}\to 0.    
\end{align}
Notice that

\begin{align*}
|\eta_{m,\rho}-E(\eta_{m,\rho}))|^q
&=\left| \sum_{j=0}^{l}(-1)^{j+1}(\xi_{m,j,\rho}-E(\xi_{m,j,\rho}))\right|^q\\
&\leq \left(\sum_{j=0}^{l}|(\xi_{m,j,\rho}-E(\xi_{m,j,\rho})|\right)^q \\
&\leq (l+1)^{q}\left(\max_{0\leq j \leq l }|(\xi_{m,j,\rho}-E(\xi_{m,j,\rho})|\right)^q\\
&= (l+1)^{q}\max_{0\leq j \leq l }\left(|(\xi_{m,j,\rho}-E(\xi_{m,j,\rho})|^q\right).
\end{align*}
Therefore,
\begin{align}\label{ii12}
\notag E|\eta_{m,\rho}-E(\eta_{m,\rho}))|^q
\notag&\leq (l+1)^{q}E\left(\max_{0\leq j \leq l }\left(|(\xi_{m,j,\rho}-E(\xi_{m,j,\rho})|^q\right)\right)\\
\notag&\leq (l+1)^{q}E\left(\sum_{j=0}^l|(\xi_{m,j,\rho}-E(\xi_{m,j,\rho})|^q\right)\\
\notag&\leq (l+1)^{q}\sum_{j=0}^lE\left(|(\xi_{m,j,\rho}-E(\xi_{m,j,\rho})|^q\right)\\
&\leq (l+1)^{q+1}\max_{0\leq j \leq l }\left(E|(\xi_{m,j,\rho}-E(\xi_{m,j,\rho})|^q\right).
\end{align}
Furthermore, by Lemma~\ref{Hyj}, for each $0\leq j\leq l$ there exists a constant $c_j$ such that
\begin{align}\label{ii1}
E|(\xi_{m,j,\rho}-E(\xi_{m,j,\rho})|^q\leq c_j E|\xi_{m,j,\rho}|^q. 
\end{align}
It follows from Equations (\ref{ii12}) and (\ref{ii1}) that there exist a constant $c$ such that
\begin{align}\label{yaps}
E|\eta_{m,\rho}-E(\eta_{m,\rho}))|^q\leq c(l+1)^{q+1} \max_{0\leq j \leq l }E|\xi_{m,j,\rho}|^q.    
\end{align}
On the other hand, by Lemma~\ref{pq}, since $\xi_{m,j,\rho}$ counts the number of $j$-faces of $R(\mathcal{P}_n,r,\rho)$ with left most point of the $j$-face in $A\cap Q_{m,n}$, we have
\begin{align}\label{yaps2}
E|\xi_{m,j,\rho}|^q\leq \prod\limits_{i=1}^j p_i^{\binom{j+1}{i+1}} E|\xi_{m,j}|^q  
\end{align}  
where $\xi_{m,j}$ is number of $j$-faces of $R(\mathcal{P}_n,r)$ with left most point in $A\cap Q_{m,n}$. It follows, by Lemma~\ref{z} and Equation~(\ref{yaps2}), that for a constant depending on $j$ we have
\begin{align}\label{yaps3}
E|\xi_{m,j,\rho}|^q\leq \prod\limits_{i=1}^j p_i^{\binom{j+1}{i+1}} E|\xi_{m,j}|^q\leq c_j\prod\limits_{i=1}^j p_i^{\binom{j+1}{i+1}}(nr^d)^ {j+1}.    
\end{align}  
Therefore, by Equations~(\ref{yaps}), (\ref{yaps2}) and (\ref{yaps3}) combined, we have that for a constant $c$
$$E|\eta_{m,\rho}-E(\eta_{m,\rho}))|^q \leq c(l+1)^{q+1}\prod\limits_{i=1}^j p_i^{\binom{j+1}{i+1}}(nr^d)^ {j+1}. $$
Thus, by Equation~(\ref{bq2}) , we have
$$b_{q,n}\leq \frac{c(l+1)^{q+1}}{r^d}\frac{\prod\limits_{i=1}^j p_i^{\binom{j+1}{i+1}}(nr^d)^ {j+1}}{ var(\chi_A^\mathcal{P}(\rho))^\frac{q}{2}}\leq c(l+1)^{q+1}\frac{n(nr^d)^ {j}}{ var(\chi_A^\mathcal{P}(\rho))^\frac{q}{2}}  $$
and since, by Lemma~\ref{Candela11}, $var(\chi_A^\mathcal{P}(\rho))\sim \mu_{0,A} n$, then 
$$b_{q,n}\leq c(l+1)^{q+1}\frac{n(nr^d)^ {j}}{  (\mu_{0,A}n)^\frac{q}{2}}=c(l+1)^{q+1}\frac{(nr^d)^ {j}}{\mu_{0,A}^{\frac{q}{2}}  n^{\frac{q}{2}-1}}.$$
Therefore, for a different constant $c'$ we have
\begin{align}\label{muj}
b_{q,n}\leq c' \frac{(nr^d)^ {j}}{ n^{\frac{q}{2}-1}}\to 0,
\end{align}
since, by hypothesis, $nr^d\to 0$ and $q\geq 3$. Thus, by Theorem~\ref{Penrose-4}, for $A\subset \mathbb{R^d}$ is open and bounded:
\begin{align}\label{poi}
\frac{\chi_A^\mathcal{P}(\rho)-E(\chi_A^\mathcal{P}(\rho))}{\sqrt{var(\chi_A^\mathcal{P}(\rho))}}\xRightarrow{D} N(0,1).    
\end{align}
\end{proof}
We now move to the case when $A$ is $\mathbb{R}^d$.

\begin{theorem}\label{Candela9xx}
Let $\Sigma$ be one of the soft random simplicial complexes $\Sigma\in\{R(\mathcal{P}_n,r,\rho), C(\mathcal{P}_n,r,\rho)\}$ 
with multiparameter $\rho=(p_1,p_2,\dots)$ and  let $\chi^P(\rho,\Sigma)$ be defined by%
$$\chi^\mathcal{P}(\rho, \Sigma)=\sum_{i=0}^{n-1}(-1)^{i+1}f_{i}^\mathcal{P}(\rho).$$
Assume that $nr^d\to 0$ and that there exists a non-negative integer $l$ such that
$$\prod\limits_{i=1}^{l+1} p_i^{\binom{l+2}{i+1}}n^{l+2}r^{d(l+1)}\to 0~~  and~~\prod\limits_{i=1}^l p_i^{\binom{l+1}{i+1}}n^{l+1}r^{dl}\to \infty.$$
Then 
$$ \frac{\chi^\mathcal{P}(\rho, \Sigma)-E(\chi^\mathcal{P}(\rho, \Sigma))}{\sqrt{var(\chi^\mathcal{P}(\rho, \Sigma))}}\xRightarrow{D} N(0,1).$$
\end{theorem}

\begin{proof}
In this proof, we will omit several steps since it is almost identical as the proof of Theorem~\ref{Candela12xx2}, and we refer the reader to this proof for the remaining details.
We will prove the theorem for $R(\mathcal{P}_n,r,\rho)$ (the proof for $C(\mathcal{P}_n,r,\rho)$ is analogous, replacing all instances of $h$ with $\mathfrak{h}$) and  using Theorem~\ref{newcosito} instead of Proposition~\ref{Penrose-1}. 

We write $\chi^\mathcal{P}(\rho,\Sigma)$ and $f^\mathcal{P}_{i}(\rho,\Sigma)$ simply as $\chi^\mathcal{P}(\rho)$ and  $f^\mathcal{P}_{i}(\rho)$, respectively, for the remainder of the proof.
Consider the open sets
$$A_K=(-K,K)^d,~~~~A^K=\mathbb{R}^d\backslash [-K,K]^d$$
and the operator $\zeta_{n,k}^\mathcal{P}$ which assigns to a subset $A\subset\mathbb{R}^d$ the random variable $\zeta_{n,k}^\mathcal{P}(A)$ defined by
$$\zeta_{n,k}^\mathcal{P}(A)=\frac{\chi_A^\mathcal{P}(\rho)-E(\chi_A^\mathcal{P}(\rho))}{\sqrt{n}}$$
which centralizes the normalization of the random variable $\chi^\mathcal{P}(\rho, \Sigma)$ by $\sqrt{n}$ around its expected value.
Given $t\in \mathbb{R}$ and $\epsilon>0$, by Equation~(\ref{der}) we have
\begin{align}\label{edi}
|P(\zeta_{n,k}^\mathcal{P}(\mathbb{R}^d)\leq t)- P(\zeta_{n,k}^\mathcal{P}(A_K)\leq t-\epsilon)|\leq P(|\zeta_{n,k}^\mathcal{P}(A_K)-t|<\epsilon) + P(|\zeta_{n,k}^\mathcal{P}(A^K)|>\epsilon).    
\end{align}
By Chebychev's inequality, Lemma~\ref{Frieze3},  and Theorem~\ref{Candela9}, the central limit theorem established for bounded sets, the right-hand part of the inequality in Equation~(\ref{edi}) is bounded above by 
\begin{align}\label{loc}
\frac{1}{\epsilon^2}var(\zeta_{n,k}^\mathcal{P}(A^K))+ P\left(\left|(Z-t)\sqrt{\frac{var(\chi_{A_K}^\mathcal{P}(\rho))}{n}}\right|<\epsilon\right)+ c \frac{(nr^d)^ {j}}{ n^{\frac{q}{2}-1}}
\end{align}
where $Z$ is a normal standard variable. The second and third expressions in Equation~(\ref{loc}) come from the central limit theorem stated in Equation~(\ref{poi}) and the error term given in Equation~(\ref{muj}).  Notice now that, as in Equation~(\ref{der3}), since $Z$ is a standard normal variable, the second term in Equation~(\ref{loc}) can be bound by
\begin{align}\label{der31}
P\left(\left|(Z-t)\sqrt{\frac{var(\chi_{A_K}^\mathcal{P}(\rho))}{n}}\right|<\epsilon\right)\leq 2\epsilon\sqrt{\frac{n}{2\pi \cdot var(\chi_{A_K}^\mathcal{P}(\rho))}}.
\end{align}
Furthermore, by Lemma~\ref{Candela11}, $$var(\chi_{A_K}^\mathcal{P}(\rho))=\mu_{0, A_K}\cdot n$$
and 
$$var(\zeta_{n,k}^\mathcal{P}(A^K))=\mu_{0, A^K}.$$ 
Thus, from Equations~(\ref{edi}), (\ref{loc}) and (\ref{der31}) combined, and Theorem~\ref{Candela9}, the central limit theorem established for bounded sets, we have 
\begin{equation}
	\label{eq:limsup}
\limsup_{n\to\infty}\left|P(\zeta_{n,k}^\mathcal{P}(\mathbb{R}^d)\leq t)- P\left(\sqrt{\frac{var(\chi_{A_K}^\mathcal{P}(\rho)))}{n}}Z\leq t-\epsilon\right)\right|\leq \frac{1}{\epsilon^2}\mu_{0,A^K}+ 2\epsilon\sqrt{\frac{1}{2\pi \mu_{0,A_K}}}.
\end{equation}
Now, if $\Phi$ is the cumulative distribution function of the normal standard variable $Z$ then
$$P\left(\sqrt{\frac{var(\chi_{A_K}^\mathcal{P}(\rho))}{n}}Z\leq t-\epsilon\right)= \Phi\left(\sqrt{\frac{n}{var(\chi_{A_K}^\mathcal{P}(\rho))}}(t-\epsilon)\right)\to \Phi\left(\sqrt{\frac{1}{\mu_{0,A_K}}}(t-\epsilon)\right).$$
Therefore, the limit in Equation (\ref{eq:limsup}) is equal to
$$\limsup_{n\to\infty}\left|P(\zeta_{n,k}^\mathcal{P}(\mathbb{R}^d)\leq t)- \Phi\left(\sqrt{\frac{1}{\mu_{0, A_K}}}(t-\epsilon)\right)\right|\leq  \frac{1}{\epsilon^2}\mu_{0,A^K}+ 2\epsilon\sqrt{\frac{1}{2\pi \mu_{0, A_K}}}$$
Recall that, by Equation~(\ref{mu0})
$$\mu_{0,A}=\int_{A}f(x)dx$$
and therefore, by the definition of $A_K$ and $A^K$,  $\lim\limits_{K\to \infty}\mu_{0, A_K}=1$ and $\lim\limits_{K\to \infty}\mu_{0, A^K}=0$. Thus, for $K$ large enough and $\epsilon'>0$, we have
$$\limsup_{n\to\infty}\left|P(\zeta_{n,k}^\mathcal{P}(\mathbb{R}^d)\leq t)- \Phi\left(\sqrt{\frac{1}{\mu_{0, A_K}}}(t-\epsilon)\right)\right|\leq \epsilon+ 2\epsilon\sqrt{\frac{1}{2\pi(1-\epsilon') }}.$$
then 
$$\limsup_{n\to\infty}\left|\left|P(\zeta_{n,k}^\mathcal{P}(\mathbb{R}^d)\leq t)-\Phi\left(t-\epsilon\right)\right|-\left| \Phi\left(\sqrt{\frac{1}{\mu_{0, A_K}}}(t-\epsilon)\right)-\Phi\left(t-\epsilon\right)\right|\right|\leq \epsilon+ 2\epsilon\sqrt{\frac{1}{2\pi(1-\epsilon') }}.$$
since $\Phi\left(\sqrt{\frac{1}{\mu_{0, A_K}}}(t-\epsilon)\right)\to \Phi\left(t-\epsilon\right)$ as $K\to\infty$, it follows that for $\epsilon_1>0$ we have
$$\limsup_{n\to\infty}\left|\left|P(\zeta_{n,k}^\mathcal{P}(\mathbb{R}^d)\leq t)- \Phi\left(t-\epsilon\right)\right|-\epsilon_1\right|\leq \epsilon+ 2\epsilon\sqrt{\frac{1}{2\pi(1-\epsilon') }}+ 2\epsilon_1.$$
Finally, since $\Phi\left(t-\epsilon\right)\to \Phi\left(t\right)$ as $\epsilon\to 0$ and $\epsilon, \epsilon'$ and $\epsilon_1$ were arbitrary, we have
$$\limsup_{n\to\infty}\left|P(\zeta_{n,k}^\mathcal{P}(\mathbb{R}^d)\leq t)- \Phi\left(t\right)\right|=0$$
and then, 
$$\lim_{n\to\infty}\left|P(\zeta_{n,k}^\mathcal{P}(\mathbb{R}^d)\leq t)- \Phi\left(t\right)\right|=0.$$
Therefore,
$$\zeta_{n,k}^\mathcal{P}\xRightarrow{D} N(0,1)$$
and this is
$$ \frac{\chi^\mathcal{P}(\rho)-E(\chi^\mathcal{P}(\rho))}{\sqrt{var(\chi^\mathcal{P}(\rho))}}\xRightarrow{D} N(0,1).$$

\end{proof}
\subsection{De-Poissonization}
Finally, the remaining work is to use a de-Poissonization argument to guarantee that the central limit theorem proved above for $\chi_A(\rho)$  is also true for $\chi(\rho)$ itself.
\begin{theorem}\label{Candela16}
Let $\Gamma$ one the soft random simplicial complexes $\Gamma\in\{R(n,r,\rho), C(n,r,\rho)\}$ with multiparameter $\rho=(p_1,p_2,...,p_{n-1},...)$ and let $\chi(\rho,\Gamma)$ be its Euler characteristic. Assume that $nr^d\to 0$ and that there exists a non-negative integer $l$ such that 
$$\prod\limits_{i=1}^{l+1} p_i^{\binom{l+2}{i+1}}n^{l+2}r^{d(l+1)}\to 0~~  and~~\prod\limits_{i=1}^l p_i^{\binom{l+1}{i+1}}n^{l+1}r^{dl}\to \infty$$
then 
$$ \frac{\chi(\rho,\Gamma)-E(\chi(\rho,\Gamma))}{\sqrt{n}}\xRightarrow{D} N(0,1).$$
\end{theorem}
\begin{theorem}\label{Candela15}
Let $R(n,r,\rho)$ be the soft random Vietoris-Rips complex and $C(n,r,\rho)$ the random soft \v{C}ech complex with multiparameter $\rho=(p_1,p_2,...,p_{n-1},...)$. Let $f_k(\rho)$ the number of $k$-faces of any of these two models. Then if $nr^d\to 0$  and $\prod\limits_{i=1}^k p_i^{\binom{k+1}{i+1}}n^{k+1}r^{dk}\to \infty$:
$$\frac{ f_k(\rho, R(n,r,\rho))-E(f_k(\rho, R(n,r,\rho)))}{\sqrt{\mu_{k+1}\prod\limits_{i=1}^kp_i^{\binom{k+1}{i+1}}n^{k+1}r^{dk}}}\xRightarrow{D} N(0, 1)$$
and
$$\frac{ f_k(\rho, C(n,r,\rho))-E(f_k(\rho, C(n,r,\rho)))}{\sqrt{\nu_{k+1}\prod\limits_{i=1}^kp_i^{\binom{k+1}{i+1}}n^{k+1}r^{dk}}}\xRightarrow{D} N(0, 1)$$
\end{theorem}
\begin{proof}[Proof of Theorem~\ref{Candela16} ]
We will prove the theorem for $R(\mathcal{P}_n,r,\rho)$ (the proof for $C(\mathcal{P}_n,r,\rho)$ is analogous, replacing all instances of $h$ with $\mathfrak{h}$) and  using Theorem~\ref{newcosito} instead of Proposition~\ref{Penrose-1}. We write $\chi(\rho,\Sigma)$ and $f_{i}(\rho,\Sigma)$ simply as $\chi(\rho)$ and  $f_{i}(\rho)$, respectively, for the remainder of the proof.

By Lemma~\ref{Candela11},  under the hypothesis of  Theorem~\ref{Candela16}, a.a.s we have
$$\chi^\mathcal{P}(\rho)=f_{0}^\mathcal{P}(\rho)-f_{1}^\mathcal{P}(\rho)+f_{2}^\mathcal{P}(\rho)-f_{3}^\mathcal{P}(\rho,)+...+(-1)^{l+1}f_{l}^\mathcal{P}(\rho)=\sum_{i=0}^l(-1)^{i+1}f_{i}^\mathcal{P}(\rho).$$
Also, since
$$f_{k}(\rho)=\sum_{Y\subset \mathcal{X}_n}h_{n,r,k+1,\rho}(Y),$$
we have that
$$\chi(\rho)=\sum_{j=0}^{l}\sum_{Y\subset \mathcal{X}_n, |Y|=j+1}(-1)^{j+1}h_{n,r,j+1,\rho}(Y),~~~a.a.s.$$
Again, we will use Theorem~\ref{Penrose-6} to recover the central limit theorem. Using the conventions from Definition \ref{h's}, for $Y\subset \mathcal{X}_n$ with $|Y|=j+1$, define  the following functional

\begin{align}\label{aguanta2}
H_n( \mathcal{X}_n)=\sqrt{n}\frac{\sum\limits_{j=0}^{l}\sum\limits_{Y\subset \mathcal{X}_n}(-1)^{j+1}h_{n,r,j+1,\rho}(Y)}{\sqrt{n}}=\sum\limits_{j=0}^{l}\sum\limits_{Y\subset \mathcal{X}_n}(-1)^{j+1}h_{n,r,j+1,\rho}(Y).    
\end{align}
Consider the increments $R_{m,n}$ of the functional $H$ given by 
\begin{align*}
R_{m,n}&=H_n(\mathcal{X}_{m+1})- H_n(\mathcal{X}_{m})\\
&=\sum\limits_{j=0}^{l}(-1)^{j+1}\left(\sum\limits_{Y\subset \mathcal{X}_{m+1}}h_{n,r,j+1,\rho}(Y)- \sum\limits_{Y\subset \mathcal{X}_m}h_{n,r,j+1,\rho}(Y)\right)\\
&=\sum\limits_{j=0}^{l}(-1)^{j+1}D_{m,n}^j    
\end{align*}
where
$$D_{m,n}^j=\sum\limits_{Y\subset \mathcal{X}_{m+1}}h_{j+1,\rho}(Y)- \sum\limits_{Y\subset \mathcal{X}_m}h_{j+1,\rho}(Y).$$
Notice that $D_{m,n}^j$ counts the number of $j$ faces in $R(\mathcal{X}_{m+1},r,\rho)$ that contain the vertex $X_{m+1}$. Thus, from Equation~(\ref{E1}), uniformly over $n-n^{\frac{2}{3}}\leq  m \leq n+n^{\frac{2}{3}}$, we have
$$E(D_{m,n}^j)\sim (j+1)\prod\limits_{i=1}^jp_i^{\binom{j+1}{i+1}}(nr^d)^{j}\mu_{j+1}$$
Therefore, by definition of $R_{m,n}$,
\begin{align*}
E(R_{m,n})\sim\sum\limits_{j=0}^{l}(-1)^{j+1} (j+1)\prod\limits_{i=1}^jp_i^{\binom{j+1}{i+1}}(nr^d)^{j}\mu_{j+1} \to 0
\end{align*}
since by hypothesis $nr^d\to 0$. Therefore, for $\alpha=0$ and $\gamma=\frac{2}{3}$, we have 
\begin{equation}\label{lim4}
\lim_{n\to\infty}\left( \sup_{n-n^\gamma\leq m\leq n+n^\gamma}|E(R_{m,n})-\alpha|)\right)\to 0.
\end{equation}
We now examine $E(R_{m_1,n}R_{m_2,n})$, for $m_1<m_2$,
\begin{align*}
E(R_{m_1,n}R_{m_2,n})&=E\left(\left(\sum\limits_{j=0}^{l}(-1)^{j+1}D_{m_1,n}^j\right)\left(\sum\limits_{s=0}^{l}(-1)^{s+1}D_{m_2,n}^s\right)\right)\\
&=E\left(\sum\limits_{j,s=0}^{l}(-1)^{s+j+2}D_{m_1,n}^jD_{m_2,n}^s\right)\\
&=\sum\limits_{j,s=0}^{l}(-1)^{s+j+2}E(D_{m_1,n}^jD_{m_2,n}^s).  
\end{align*}
On the other hand, using the same calculation used to deduce (\ref{proD's}), we can obtain that, uniformly over $m_1,m_2\in [n-n^{\frac{2}{3}}, n+n^{\frac{2}{3}}]$, as $n\to \infty$
$$E(D_{m,n}^jD_{m',n}^s)\sim(j+1)(s+1)\mu_{j+1}\mu_{s+1}(nr^d)^{j+s}\prod\limits_{i=1}^{\max(s,j)}p_i^{\binom{s+1}{i+1}+ \binom{j+1}{i+1}}.$$
Thus
$$E(R_{m_1,n}R_{m_2,n})\sim\sum\limits_{j,s=0}^{l}(-1)^{s+j+2} (j+1)(s+1)\mu_{j+1}\mu_{s+1}(nr^d)^{j+s}\prod\limits_{i=1}^{\max(s,j)}p_i^{\binom{s+1}{i+1}+ \binom{j+1}{i+1}}$$
and since $nr^d\to 0$, it follows that
$$ E(R_{m,n}R_{m,n})\to 0.$$
Therefore, for $\alpha=0$ and $\gamma=\frac{2}{3}$, we have
\begin{equation}\label{lim5}
\lim_{n\to\infty}\left( \sup_{n-n^\gamma\leq m_1<m_2\leq n+n^\gamma}|E(R_{m_1,n}E(R_{m_2,n})-\alpha^2|)\right)\to 0.
\end{equation}
Furthermore, uniformly over $m\in [n-n^{\frac{2}{3}}, n+n^{\frac{2}{3}}]$, as $n\to \infty$

\begin{align*}
\frac{1}{\sqrt{n}}E(R_{m,n}^2) &\sim \sum\limits_{j,s=0}^{l}(-1)^{s+j+2} (j+1)(s+1)\mu_{j+1}\mu_{s+1}\frac{(nr^d)^{j+s}}{\sqrt{n}}\prod\limits_{i=1}^{\max(s,j)}p_i^{\binom{s+1}{i+1}+ \binom{j+1}{i+1}}\\
&\sim \sum\limits_{j,s=0}^{l}(-1)^{s+j+2} (j+1)(s+1)\mu_{j+1}\mu_{s+1}n^{j+s-\frac{1}{2}}r^{d(j+s)}\prod\limits_{i=1}^{\max(s,j)}p_i^{\binom{s+1}{i+1}+ \binom{j+1}{i+1}}\\
&\leq \sum\limits_{j,s=0}^{l}(-1)^{s+j+2} (j+1)(s+1)\mu_{j+1}\mu_{s+1}n^{j+s-\frac{1}{2}}r^{d(j+s-\frac{1}{2})}\prod\limits_{i=1}^{\max(s,j)}p_i^{\binom{s+1}{i+1}+ \binom{j+1}{i+1}}\to 0 
\end{align*}
Therefore, 
\begin{equation}\label{lim6}
\lim_{n\to\infty}\left( \frac{1}{\sqrt{n}}\sup_{n-n^\gamma\leq m\leq n+n^\gamma}|E(R_{m,n}^2|)\right)\to 0.
\end{equation}
By estimates (\ref{lim4}), (\ref{lim5}) and (\ref{lim6}), together with Theorem~\ref{Candela9xx}, the functional $H$ defined in Equation~(\ref{aguanta2}) satisfies all the conditions of the Theorem~\ref{Penrose-6} for $\alpha=0$. Therefore, by Theorem~\ref{Penrose-6}, for $\sigma$ such that
$$\frac{var(H_n( \mathcal{P}_n))}{n}\to \sigma^2$$
we have
$$\frac{H_n( \mathcal{X}_n)- E(H_n( \mathcal{X}_n))}{\sqrt{n}}\Rightarrow N(0, \sigma^2).$$
That is
$$\frac{ \chi(\rho)- E(\chi(\rho))}{\sqrt{n}}\xRightarrow{D} N(0, \sigma^2).$$
Finally, for Theorem~\ref{Candela9xx}, we have
\begin{align*}
\frac{var(H_n( \mathcal{P}_n))}{n}&=\frac{1}{n} var\left(\frac{\sqrt{n}(\chi^{\mathcal{P}}(\rho)-E(\chi^{\mathcal{P}}(\rho))) }{\sqrt{n}}\right)\\
&= var\left(\frac{\chi^{\mathcal{P}}(\rho)-E(\chi^{\mathcal{P}}(\rho)) }{\sqrt{n}}\right)\\
&=1. 
\end{align*}
Therefore,
$$\frac{\chi(\rho)-E(\chi(\rho)) }{var(\chi(\rho))}\xRightarrow{D} N(0, 1),$$
as desired. 
\end{proof}

\section*{Acknowledgements}
I want to thank my advisor Antonio Rieser for guiding me throughout my research and for sharing his knowledge with me. Additionally, I would like to thank   Octavio Arizmendi and Arturo Jaramillo for helpful comments on earlier versions of the manuscript.

\appendix

\section{Inequalities}

\begin{lemma}[Lemma 20.1 \cite{MR3675279}]\label{Frieze1} Let $X$ be  a non-negative random variable. Then, for all $t>0$:
$$P(X\geq t)\leq \frac{E(X)}{t}.$$
\end{lemma}

\begin{lemma}[Lemma 20.2 \cite{MR3675279}]\label{Frieze2} Let $X$ be  a non-negative integer valued random variable. Then:
$$P(X>0)\leq E(X).$$
\end{lemma}

\begin{lemma}[Lemma 20.3 \cite{MR3675279}]\label{Frieze3} If $X$ is a random variable with a finite mean and variance. Then, for all $t>0$:
$$P(|X-E(X)|\geq t)\leq \frac{var(X)}{t^2}.$$
\end{lemma}

\begin{lemma}[Lemma 20.4 \cite{MR3675279}]\label{Frieze4} If X is a non-negative integer valued random variable then
$$P(X=0)\leq \frac{var(X)}{(E(X))^2}=\frac{E(X^2)}{(E(X))^2}-1.$$
\end{lemma}

\begin{lemma}[Theorem 1.6.2 \cite{MR2722836}]\label{Jensen} Suppose $\phi$ es convex function. Then 
$$E(\phi(X))\geq \phi(E(X))$$
provided both expectations exist i.e $E|X|$ and $E|\phi(X)|<\infty$. 
\end{lemma}

\begin{lemma}[Theorem 1.6.3 \cite{MR2722836}]\label{Holder} If $p,q\in[1,\infty]$ with $\frac{1}{p}+\frac{1}{q}=1$, then
$$E|XY|\leq E(|X|^p)^{\frac{1}{p}}E(|Y|^q)^{\frac{1}{q}}.$$
\end{lemma}

\end{document}